\input amstex
\documentstyle{amsppt}
\magnification=\magstep1
\vsize=19cm

\topmatter

\title{On locally trivial extensions of topological spaces by a pseudogroup} 
\endtitle

\rightheadtext{On extensions of spaces by a pseudogroup}

\author  Andre Haefliger and Ana Maria Porto F. Silva\endauthor

\abstract{In this paper we restrict ourselves to the particular case where the pseudogroup is   $\Gamma \ltimes G$  given  by the action of a dense subgroup $\Gamma$ of a Lie group $G$ acting on $G$ by left translations. 
For a Riemannian foliation $\Cal F$ on a complete Riemannian manifold $M$ which is transversally parallelizable in the sense of Molino, let $X$ be the space of leaves closures. The holonomy pseudogroup of $\Cal F$ is an example of a locally trivial extension of $X$ by $\Gamma \ltimes G$. The study of a generalization of this particular case shall be the main purpose of this paper. }

 \endabstract

\endtopmatter

\heading{Introduction 
}\endheading

Inspired by Molino theory of Riemmanian foliations [6], the general study of pseudogroups of local isometries was sketched in [4]. In the present paper we restrict ourselves to the study of the holonomy pseudogroups which appear in the particular case of Riemannian foliations $\Cal F$ which are transversally complete (see [5]). In this case, the holonomy pseudogroup of the closure of a leaf is the pseudogroup $\Gamma \ltimes G$ given by the action of a dense subgroup $\Gamma$ of a Lie group $G$ acting on $G$ by left translations. Let $X$ be the space of orbit closures of the holonomy pseudogroup $\Cal G$ of $\Cal F$. Then, in our terminology, $\Cal G$ is a locally trivial extension of $X$ by $\Gamma \ltimes  G$ and it is called a $(\Gamma \ltimes G)$-extension of $X$.
We write $\rho : \Gamma \to G$ for the inclusion of $\Gamma$ in $G$.

This paper presents most of  the results contained  in the thesis of the second author [10] and announced in [11], but, sometimes, with different notations and other proofs.

Roughly speaking, in our setting, a $(\Gamma \ltimes G)$-extension of a  space $X$ is a toplogical groupoid $\Cal G$ with a continuous projection $p: \Cal G \to X$ which can be considered as a fibration with fiber $\Gamma \ltimes G$ associated to an $Aut (\Gamma \ltimes G)$-principal bundle, where $Aut(\Gamma \ltimes G)$ is the group of automorphisms of the groupoid $\Gamma \ltimes G$.

The group $Aut(\Gamma \ltimes G)$ is a Lie group isomorphic to $Aut(\rho) \ltimes G$, where $G$ is with its usual topology and $Aut(\rho)$ is the discrete group whose elements are the automorphisms $\overline{\alpha}$ of $G$ such that $\overline \alpha(\rho(\gamma)) g = \rho(\alpha (\gamma))\overline \alpha (g)$.

 The map $Aut(\rho) \ltimes G \to Aut(\Gamma \ltimes G)$
 sending $(\alpha, g) \in Aut(\rho) \ltimes G$ to the automorphism 
$(\gamma,g') \to (\alpha(\gamma),\overline \alpha(g')g'g^{-1})$ is an isomorphism.

The computation of $Aut(\rho)$ often involves  number theory considerations depending of the nature of $G$, compact, nilpotent,  abelian or solvable. Many examples of the determination of this group shall appear throughout this paper, in particular in part III. 

The quotient of $Aut(\rho)$ by the subgroup $Int(\rho)$ of its inner automorphisms is noted $Out(\rho)$.

The main goal of our work is the construction and the determination of  the homotopy type of a classifying space for $(\Gamma \ltimes G)$- extensions of topological spaces.

As a sample of such results, we state  the following theorem (see in part II, theorem 4.1, for a more precise statement):

\proclaim {Theorem} The space $BAut(\Gamma \ltimes G)$, also noted $X^{\Gamma \ltimes G}$,  is a classifying space for the set of homotopy classes of $(\Gamma \ltimes G)$-extensions of paracompact spaces.
Moreover the corresponding $(\Gamma \ltimes G)$-extension $p^{\Gamma \ltimes G}:B\Cal G^{\Gamma \ltimes G}\to 
X^{\Gamma \ltimes G}$ is a fiber bundle with fiber $B(\Gamma \ltimes  G)$.

Also the space $ B\Cal G^{\Gamma \ltimes G}$ has the homotopy type an Eilenberg-Maclane complex $K(Aut(\rho),1)$.
\endproclaim

We now give a brief description of the content of the paper

In  part I, called Basic Notions on Groupoids, we recall the  definition of pseudogroups and the equivalent notion of their groupoid of germs. We define the notion of weak isomorphisms among topological groupoids (called equivalences in many papers like [2] and [4] ). We give many examples of such isomorphisms. We also recall the construction of the  holonomy groupoid of a foliation $\Cal F$ on a manifold $M$ using a submanifold $T$ of $M$ transversal to the leaves and cutting all of them. The holonomy groupoids using two such transversals are weakly isomorphic.

The core of the paper is contained in part II, entitled: $(\Gamma \ltimes G)$-extensions of topological spaces. Note that for us the group $\Gamma$ can be any dense subgroup of $G$. For instance, it could be the group $G$ with the discrete topology. In the case of Riemannian foliations, the group $\Gamma$ is always countable.

 In section 1 of part II, we  give a  precise description of the group $Aut (\Gamma \ltimes G)$ of automorphisms of $\Gamma \ltimes G$ as mentionned above, and many examples of the computation of the group $Aut(\rho)$ in 1.11.

An important homomorphism $\mu: \Gamma \to Aut(\Gamma \ltimes G) = Aut(\rho) \ltimes G$ is defined 
by $\mu (\gamma) = (ad(\gamma), \rho(\gamma)^{-1})$. The subgroup $\mu(\Gamma)$ is a normal subgroup, but in general, it is not a closed subroup. Therefore,
the quotient $Aut(\Gamma \ltimes G) / \mu(\Gamma)$ noted $Out(\Gamma \ltimes G)$ is not in general a Lie group.

Section 2 of Part II contains the basic definitions and properties of locally trivial $(\Gamma \ltimes G)$-extensions of topological spaces. Such an extension $p:\Cal G \to X$ can be described as follows.
Each point $x \in X$ has a neighbourhood $U$ such that there is a weak isomorphism $E_U$ from $p^{-1} (U)$  to the groupoid $\Gamma \ltimes (G \times U)$ given by the action of $\gamma \in\Gamma$ on $G \times U$ defined by $(g,u) \mapsto (\rho(\gamma)g,u)$. Moreover, if $U'$ is an other such open set containing $x$, then
$E_{U'} \circ {E_U}^{-1}$ is given by a continuous map $\phi_{U'U}: U \cap U' \to Aut(\Gamma \ltimes G)$.

Two such $(\Gamma \ltimes G)$-extensions  $p_1: \Cal G_1 \to X$ and $p_2: \Cal G_2 \to X$ are weakly isomorphic if there is a weak isomorphism from $\Cal G_1$ to $\Cal G_2$ commuting with $p_1$ and $p_2$.

Let $(U_i)_{i \in I}$ be a fine enough open cover such that $E_{U_i} \circ {E_{U_j}}^{-1}$ is given by a continuous map $\phi_{ij}:U_i \cap U_j  \to Aut(\Gamma \ltimes G)$. Then, $\phi_{ij}$ is a $1$-cocycle with values in the sheaf $\tilde{Aut}(\Gamma \ltimes G)$ of germs of continuous maps of open sets of $X$ to the group $Aut(\Gamma \ltimes G)$. We prove that the set of weak isomorphisms classes of  $(\Gamma \ltimes G)$-extensions of $X$ is in bijection with the set $H^1(X, \tilde {Aut}(\Gamma \ltimes G))$ (theorem 2.5).

 The notion of equivalence between  $(\Gamma \ltimes G)$-extensions is also defined in 2.1.
In theorem 2.7, we prove that   the
   set of equivalence classes of $(\Gamma \ltimes G)$-extensions of $X$ is equal to the set $H^1(X, \tilde {Out}(\Gamma \ltimes G))$. There, $\tilde {Out}(\Gamma \ltimes G)$ is defined as the quotient of the sheaf $\tilde {Aut} (\Gamma \ltimes G)$ by the subsheaf $\mu (\tilde \Gamma)$.

  At the end of section 2 we define the notion of homotopy  for $(\Gamma \ltimes G)$-extensions of $X$ preserving a base point $* \in X$. For a sheaf $\Cal A$ over $X$ we define $\Cal A_*$ as the subsheaf of $\Cal A$ such that the stalk of $\Cal A_*$ over $*$ is just the unit element and prove that the set of homotopy classes of $(\Gamma \ltimes G)$-extensions preserving the base point $*$ is equal to $H^1(X, \tilde{Aut}(\Gamma \ltimes G)_*)$. Similarly,  $H^1(X, \tilde {Out}(\Gamma \ltimes G)$ is the set of homotopy classes of pointed equivalence classes of equivalences of  $(\Gamma \ltimes G)$ (see theorem 2.11).

  In section 3, we recall the existence of a classifying space for a topological groupoid $G$ extending  Milnor's join construction for a topological group (see [1] and [3]). Given a topological space $X$, then $\tilde G$ is as above the sheaf over $X$ of germs of continuous maps from open sets of $X$ to $G$. Then, we generalize these notions to take care of base points; a base point   in $G$ is determined  by the choice of a point $*$ in the space of objects of $G$. Given a base point $*$ in a space, then $\tilde G_*$ is the subsheaf of $\tilde G$ whose stalk above $*$ contains only the germ of the constant map from $X$ to $*$. We define $[X,BG]_*$ as the set of homotopy classes of continuous maps from $X$ to $BG$ preserving the base points and we get  the  theorem 3.3, namely $H^1(X,\tilde G_*) = [X, BG]_*$. 
  
  In 3.7, a fondamental lemma  is proved; this is a crucial step used in the proof of theorem 4.1.

In the last section 4 of part II, we prove the main results of this paper.  
  First, we prove in theorem 4.1  the existence of a universal pointed  $(\Gamma \ltimes G)$-extension:  $p^{\Gamma \ltimes G}: \Cal G^{\Gamma \ltimes G} \to X^{\Gamma \ltimes G}$ such that any pointed $(\Gamma \ltimes G)$-extension $p: \Cal G \to X$ is (weakly) isomorphic to the pulback of the  universal extension by a base point preserving map $f: X \to X^{\Gamma \ltimes G}$. We then show that $B\Cal G^{\Gamma \ltimes G}$ has the homotopy type of the Eilenberg-MacLane complex $K(Aut(\rho),1)$. 
  
  As a corollary, we get the determination of all the homotopy groups $\pi_i(X^{\Gamma \ltimes G})$ of the classifying space. Namely, they are equal to $Out(\rho)$  for $i=1$, to $\Gamma_0$ for $i =2$, to $\pi_i(BG,*)$ for $i>2$.

 In theorem 4.4 we describe the Postnikov  decomposition of $X^{\Gamma \ltimes G}$.

  In part III, we illustrate those results  by constructing several examples. In section 1, we consider the case   where the center $\Gamma_0$ of $\Gamma$ is a discrete subgroup of $G$. This always the case if $G$ is a compact group.  We prove that $\mu(\Gamma)$ is a discrete subgroup of $G$. Therefore, $Out(\Gamma \ltimes G)$ is a Lie group. We prove  that this group is isomorphic to $Out(\rho) \ltimes G/\rho(\Gamma _0)$.

In section 2, we assume that  $G$ is $R$ or $C$ and that $\Gamma$ is of finite rank. In that case, the theory of number fields plays a central role.

In section 3, we consider the case where $\rho(\Gamma)$ is a finitely generated subgroup of a nilpotent Lie group $G$. As $G$ is assumed to be simply connected, the group $\Gamma$ is without torsion. We can, then, apply  Malcev theory to obtain several interesting results.

In section 4, we consider the case where $G$ is solvable.

At the end of the paper,  we give in a postface indications for many possible generalizations of the results contained in this paper.

In our paper, we have restricted ourselves strictely to the case mentionned in the abstract.
This paper is mostly self-contained. The only references used  to prove  the general results   are Serre [9] for sheaf theory, [1], [3] and [4]. 
Other references are given at appropriate places for some applications or examples.

For an introduction to the theory of pseudogroup of local isometries and its relation with Molino theory which is the basic ground of this paper, see the excellent survey [8] of E. Salem. The notion of homotopy groups of pseudogroups is also defined in this paper.

\heading{\bf TABLE OF CONTENTS.}\endheading

I. Basic notions on groupoids.

II. ($\Gamma \ltimes G$)-extensions of topological spaces.

III. Examples

Postface.

Bibliography.

\newpage

\heading{I. Basic notions on groupoids}\endheading

\subheading{1. Pseudogroups and groupoids}

In this section we introduce the notion of pseudogroups and the notion of its associated groupoid of germs.

\noindent {\bf 1.1. Pseudgroups of local homeomorphisms}
 
 A {\it local homeomorphism}  of a topological space $X$ is a homeomorphism $h: U \to V$, where $U$ and $V$ are open sets of $X$. The open set $U$ is called the {\it source} of $h$ and $V$ the {\it target} of $h$. The composition of two such homeomorphisms $h_1:U_1 \to V_1$ and $h_2:U_2 \to V_2$ (noted $h_2 \circ h_1$) is defined if $V_1 \cap U_2 \neq \emptyset$ and is the local homeomorphism $U_1 \cap h_1^{-1}(V_1\cap U_2) \to h_2(V_1 \cap U_2)$ mapping $x$ to $h_2(h_1(x)))$.

 A {\it pseudogroup} $P$ of local homeomorphisms of $X$ is a collection of local homeomorphisms of $X$ such that :
 
 1) The identity homeomorphism of $X$ belongs to $P$.
 
 2) If $h$ belongs to $P$, then the inverse $h^{-1}$ belongs to $P$ as well as the restriction of $h$ to any open set contained in its source.
 
 3) If $h_1$ and $h_2$ belong to $P$, then $h_1 \circ h_2$ belongs to $P$ if it is defined.
 
  Two points $x$ and $y$ in $X$ are in the same {\it orbit} under $P$ if there is an $h$ in $P$ such that $y =h(x)$. The relation for two points to be in the same orbit under $P$ is an equivalence relation, in other words the orbits under $P$ form a partition of $X$. The orbit of $x$ under $P$ is noted $P.x$. The space of orbit is endowed with the quotient topology and is noted $P\backslash X$.
  
  The pseudogroup generated by a family $(f_i)_{i \in I}$ of local homeomorphisms of $X$ is the intersection of the pseudogroups containing all the elements $f_i$ of this family. Note that the pseudogroup consisting of all the local homeomorphisms of $X$ contains any pseudogroup of local homeomorphisms of $X$. As an exercise, try to describe the pseudogroup generated by two local homeomorphisms of $X$. 
  
\noindent {\bf 1.2. Groupoids}
  
  A {\it groupoid} is a small category $\Cal G$ whose all morphisms (or arrows) are invertible. Let $X$ be its set of objects. We think of an element $g \in\Cal G$ as an arrow from right to left with source $s(g) \in X$ and target $t(g) \in X$. We have an inclusion $i: X \to \Cal G$ associating to $x$ a unit $1_x $. So for an element $g \in \Cal  G$, then $1_{s(g)}$ is the right unit of $g$ and $1_{t(g)}$ is the left unit of $g$. The composition of $g_1$ and $g_2$ is defined only if $t(g_1) = s(g_2)$ and is noted $g_2g_1$. The set of composable elements is noted $\Cal G^{(2)}$. More generally for an integral number $n >0$, then $\Cal G^{(n)}$ is the subset of $\Cal G^n$ formed by  sequences $(g_n,...,g_1)$ of $n$ composable elements, namely $t(g_i) = s(g_{i+1})$ for $n$$>$ $ i \geq �1$. The composition of these elements is noted $g_n \dots g_1$.
  
  A {\it topological groupoid} is a groupoid $\Cal G$ such that $\Cal G$ and $X$ are topological spaces and  all the structure maps $s :\Cal G \to X$, $t:\Cal G \to X$, $i: X \to \Cal G$, $g \to g^{-1}$,  $(g_2,g_1) \in \Cal G^{(2)} \mapsto  g_2g_1 \in \Cal G$ are continuous. The topology of $\Cal G^{(2)}$ is induced by the topology of $\Cal G^2$ and the topology of $X$ is induced by the inclusion $i: X \to \Cal G$.
  
 Let $\Cal G$ and $\Cal G'$ be two topological groupoids with spaces of objects $X$ and $X'$ respectively. A {\it homomorphism} $\phi: \Cal G \to \Cal G'$ is a continuous functor from $\Cal G$ to $\Cal G'$. It induces a continuous map $\phi_0:X \to X'$ defined by $\phi(1_x) = 1_{\phi_0(x)}$.
The set of homomorphisms from $\Cal G$ to $\Cal G'$ is noted $Hom(\Cal G,\Cal G')$.
 If $\phi'$ is a  homomorphism from $\Cal G'$ to a  groupoid $\Cal G''$, then the composition $\phi'\phi$ is defined as the composition of functors, i.e $(\phi'\phi)(g) = \phi'(\phi(g)),\ \forall g \in \Cal G$. 
  
  An {\it \'etale groupoid} is a topological groupoid $\Cal G$ such that the maps $s,t : \Cal G \to X$ are \'etale, i.e.  each $g \in \Cal  G$ is contained in an open set such that the restrictions of $s$ and $t$ to this set is a homeomorphism onto an open set of $X$.
  
  As an example of an \'etale groupoid consider a group $\Gamma$ acting by homeomorphisms on a topological space $X$. The image of $x\in X$ by $\gamma\in \Gamma$ is noted $\gamma.x$. We can associate to this action an \'etale groupoid $\Gamma \ltimes X$ homeomorphic to $\Gamma \times X$ (note that $\Gamma$ has the discrete topology). The source and target projections map $(\gamma,x)$ to $x$ and $\gamma.x$ respectively. The composition is defined by $(\gamma_2,\gamma_1.x)(\gamma_1,x)= (\gamma_2\gamma_1,x)$. The inclusion $X \to \Gamma \ltimes X$ maps $x$ to $(1,x)$.

 \noindent {\bf 1.3. The space of germs}
  
 Let $X$ and $Y$ be two topological spaces and let  $f: U \to Y$ and $f': U' \to Y$ be two continuous maps, where $U$ and $U'$ are open sets of $X$. They have the same {\it germ} at a point $x \in U \cap U'$   if their restrictions to an open nbhd of $x$ are equal. The relation for two such maps to have the same germ at a point $x$ of their source is an equivalence relation. The equivalence class of $f$ is called the germ of $f$ at $x$ and is noted $f^x$. The point $x$ is called the source $s(f^x)$ of $f^x $and the point $t(f^x) = f(x)$ its target. More generally for a continuous map $f:U \to Y$ we note $f^{U}$ the set of germs of $f$ at the various points of $U$. The map  associating to $x\in U$ the germ of $f$ at $x$ is a bijection from
 $U$ onto $f^{U}$ .

 Let $J^0(Y,X)$ be the set of germs of all continuous maps from $X$ to $Y$ at the various points of $X$. The natural topology on this set is the topology whose open sets are the sets $f^{U}$. The map $s:J^0(Y,X) \to X$ is an \'etale map, and the map $t: J^0(X,Y)\to Y$ is continuous. 
 
 Let $Z$ be a topological space. Let $J^0(Z,Y) \times_Y J^0(Y,X)$ be the subspace of $J^0(Z,Y) \times J^0(Y,X)$ formed by pairs $(g^y,f^x)$ such that $y =f(x)$, where $g$ is a continuous map from an open set of $Y$ to $Z$. The map $J^0(Z,Y) \times_Y J^0(Y,X) \to J^0(Z,X)$ associating to a pair $(g^y,f^x)$ the germ of $ g \circ f$ at $x$ is continuous.
 
 Consider the open subspace $\tilde J^0(X)$ of $J^0(X,X)$ formed by the germs of local homeomorphisms of $X$. With respect to the composition of germs, $\tilde J^0(X)$ is an \'etale groupoid.
 
 To each \'etale groupoid $\Cal G$  with space of objects $X$, we can associate an \'etale groupoid $\tilde \Cal G$ with space of objects $X$ which is a subgroupoid of $\tilde J^0(X)$,
  and a surjective homomorphism $j: \Cal G \to \tilde \Cal G$ projecting to the identity of $X$ defined as follows. 
  
  Given $g \in \Cal G$ with source $x$, one can find a local section $\sigma$ of the source map such that $\sigma(x) =g$. Then $j(g)$ is the germ of $\sigma$ at $x$. It is clear that $j$ is a continuous homomorphism of \'etale groupoid.
   The kernel of $j$ is called the inertia subgroupoid of $\Cal G$ ; its elements are the elements $g \in \Cal G$ such that $s(g)=t(g)$. 
 
 As an example consider the case where $\Cal G = \Gamma \ltimes X$ (see 1.2). Then $\tilde \Cal G$ is the identity of $X$ if and only if $\Gamma$ acts simply transitively on $X$ (for instance if $X$ is a point). 
 
 As an exercise give examples where the homomorphism $j$ is  an isomorphism.
 Consider also the case where $X$ is a discrete space and where $\Gamma$ acts by permutations on $X$.

\noindent{\bf 1.4. The \'etale groupoid associated to a pseudogroup}
 
Let $P$ be a pseudogroup of local homeomorphisms of a topological space $X$. The  groupoid $\tilde P$ associated to $P$ is the \'etale groupoid whose elements are the germs of   local homeomorphisms containded in  $P$ at the various points of their source. The source and target $s,t: \tilde P \to X$ map the germ $h^x$ of an  element $h: U \to V$ of $P$  to $x$ and $h(x)$ respectively. For $h_1: U_1 \to V_1$ and $h_2: U_2 \to V_2$ and $x \in U_1 \cap h_1^{-1}(U_2)$, the composition ${h_2}^{h_1(x)}{h_1}^x$ is defined as ${(h_2\circ h_1)}^x$. The inclusion $i: X \to \tilde P$ maps $x \in X$ to the germ at $x$ of the identity map of $X$. The \'etale groupoid associated as above to $\tilde P$ is equal to $\tilde P$.

From $\tilde P$ we can reconstruct $ P$. Indeed each element  of $\tilde P$ with source $x$ is the germ at $x$ of an \'etale map $f : u \mapsto t(\sigma(u))$, where $\sigma : U \to \tilde P$ is defined as in 1.3 on an open nbhd $U$ of $x$. The restriction of $f$ to a small enough neighbourhood $U'$ of $x$ is a local homeomorphism $f': U' \to V'$. Then $P$ is generated by all such local homeomorphisms. For instance the pseudogroup associated to the \'etale groupoid $\tilde J^0(X)$ is the pseudogroup of all local homeomorphisms of $X$.

A basic example of pseudogroup which will be important for us is the following. Let $\Gamma$ be a discrete group acting by homeomorphisms on a connected topological space $X$. Let us assume  that  the action is {\it quasi-analytic} ; this means that if $\gamma \in \Gamma$ restricted to an open set is the identity, then $\gamma$ is the identity element of $\Gamma$. Let $P$ be the pseudogroup of local homeomorphisms of $X$ generated by the elements of $\Gamma$. We note $\tau_\gamma$ the homeomorphism of $X$ defined by $x \mapsto \gamma.x$.
Then 
\proclaim{Lemma}The homomorphism of etale  groupoids $\phi: \Gamma \ltimes X \to \tilde P$ sending $(\gamma,x)$ to the germ $\tau^x_\gamma$ of $\tau_\gamma$ at $x$ is an isomorphism.

  \endproclaim
\demo{Proof}  The map $\phi$  is continuous by definition of the germ topology. It is  surjective because $P$ is generated by the action of the group $\Gamma$. It is also injective ; indeed if
 $\tau_\gamma^x = \tau_{\gamma'}^x $ for some $x$, then there is an open nbhd $U$ of $x$ such that $\gamma.u =\gamma'.u$ for all $u \in U$. This implies $u = (\gamma^{-1}\gamma').u$ for all $u \in U$, i.e. $\gamma^{-1}\gamma'$ is the identity on $U$; as the action is quasi-analytic, we have $\gamma=\gamma'$. Therefore $\phi$ is an isomorphism.
  $\square$
\enddemo

\noindent {\bf 1.5. Weak isomorphisms of pseudogroups}

Let $P_1$ and $P_2$ be pseudgroups of local homeomorphisms of spaces $X_1$ and $X_2$ respectively. A {\it weak isomorphism} from $P_1$ to $P_2$ is a collection $E$ of homeomorphisms from open sets of $X_1$ to open sets of $X_2$ such that 

(1) for each $h_1 \in P_1$, $h_2 \in P_2$, $f \in E$, then $h_2 \circ f \circ h_1 \in E$

(2) the elements of the form   $f^{-1}\circ f'$ and $f' \circ f^{-1}$   for all $f,f' \in E$ generate $P_1$ and $P_2$ respectively.

This implies that if $ { f'}^{-1}\circ f $ is the identity of an open set $U$, then $f=f'$; otherwise $f \circ { f'}^{-1}$ would not be the identity of $f(U)$, and this would contradict (2).

A weak isomorphism  $E$ from $P_1$ to $P_2$
is noted $E: P_1 \to P_2$.

The collection $E^{-1}$ of all the inverses of elements of $E$ is a weak isomorphism from $P_2$ to $P_1$. If $E_1$ is a weak isomorphism from $P_1$ to $P_2$ and $E_2$
    is a weak isomorphism from $P_2$ to a pseudogroup $P_3$ of local homeomorphisms of $X_3$, then the collection of all compositions $f_2 \circ f_1$ of elements $f_1$ of $E_1$ with those $f_2$ of $E_2$ is a weak isomorphism from $P_1$ to $P_3$, noted $E_2 \circ E_1$. The set of weak isomorphisms from $P_1$ to $P_1$ form a group under composition of weak isomorphisms called the group of  weak automorphisms  of $P_1$.
     
     If $P_1$ and $P_2$ are pseudogroups of local diffeomorphisms of differentiable manifolds $X_1$and $X_2$, a weak isomorphism $E$ from $P_1$ to $P_2$ is called differentiable if all its elements are diffeomorphisms.
     
     A weak isomorphism is called an equivalence in [2], [3], [4], [8] and [11]. Later on from II.2.1, a notion
of equivalence  shall be defined in another context.   

 Let $\tilde E$ be the space of germs of the elements of $E$. Let $s : \tilde E \to X_1$ (resp. $t:\tilde E \to X_2$) be the map associating  to a germ at $x$ its source  (resp. its target). The maps $s$ and $t$ are \'etale maps. Let $ e \in \tilde E$ with $s(e) = x_1,\ t(e) = x_2$ and $g_1 \in \tilde P_1$ with $t(g_1) =x_1$ and $g_2 \in \tilde P_2$ with $s(g_2) =x_2$, then $g_2eg_1 $ is the element of $\tilde E$ which is the composition of germs. Note that if $g_2eg_1 = g'_2eg'_1$, then $g_2 =g'_2$ and $g_1 =g'_1$. 
 
 We have a category  whose objects are pseudogroups and whose arrows are weak isomorphisms of pseudogroups.

 In the particular case where $E: P_1 \to P_2$ is generated by a continuous map $f: X_1 \to X_2$, equivalently if there is a continuous global section $\sigma: X_1 \to \tilde E$ of the projection $s : \tilde E \to X_1$, then we get a homomorphism  $\tilde P_1 \to \tilde P_2$ defined as follows. To the germ $g_1 \in \tilde P_1$ at $u \in U$ of the  element $h_1: U \to U'$ of $ P_1$ we associate the germ  of  $f \circ h_1 \circ f^{-1}$ at $f(u)$.

Let us give a few examples of weak isomorphisms of this kind. We shall see later that we can always be in the case mentionned above after replacing  $\tilde P_1$ by a groupoid obtained by localization of $\tilde P_1$ over a fine enough open cover of $X_1$.

(a) Let $P$ be a pseudogroup of local homeomorphisms of the topological space $X$. Let $X_0$ be an open subset of $X$ containing at least a point of each orbit under $P$. Then the collection of local homeomorphisms of $X_0$ contained in $P$ is a pseudgroup  $P_0$ weakly isomorphic to $P$. A weak isomorphism $E$ from $P_0$ to $P$ is the collection of all  elements of $ P$ with sources contained in $X_0$. The  space
of germs $\tilde E$ associated to $E$ is the subspace of the groupoid $\tilde P$ formed by the elements with sources in $X_0$.

(b) Let $\Gamma$ be a group of homeomorphisms of a topological space $X$ acting  quasi-analytically. Assume that each point $x \in X$ has an open neighbourhood $U_x$ such that $U_x \cap \gamma.U_x = \emptyset$ for all $\gamma \in \Gamma$, except for $\gamma = 1$. Then the pseudogroup $P_1$ generated by the elements of $\Gamma$ is weakly isomorphic to the pseudogroup $P_2$ generated  by the identity homeomorphism of the quotient space $\Gamma\backslash X$. A weak isomorphism $E: P_1 \to P_2$ is generated by the canonical projection $p: X \to \Gamma \backslash X$. The corresponding weak isomorphism $\tilde E: \tilde P_1 \to \tilde P_2$ is the space of germs of the projection $p$ at the various points of $X$.

(c) More generally let $f : X \to Y$ be an \'etale surjective map. Suppose that $\Gamma$ acts quasi-analytically on $X$ and leaves invariant the projection $f$. This means that, for each $x \in X$ and $\gamma \in \Gamma$, we have $\gamma.f^{-1}(f(x)) = f^{-1}f(\gamma.x)$. This gives an action of $\Gamma$ by homeomorphisms on $Y$. Let $\Lambda$ be the subgroup fixing all points of $Y$. This is a normal subgroup of $\Gamma$; the group $\Gamma/\Lambda$ acts qusi-analytically on $Y$.

 The  pseudogroup $P$ generated by the action of $\Gamma$ on $X$ is weakly isomorphic to the pseudogroup $Q$ generated by the action of $\Gamma/\Lambda$ on $Y$. A weak isomorphism $E: P \to Q$ is generated by $f$, more precisely by the restrictions of $f$ to  small enough open sets of $X$. The corresponding weak isomorphism $\tilde E$ is the space of germs of $f$ at the various points of $X$.

As an example, consider the \'etale map $R^n \to R^n/Z^n=T^n$. Let $\Gamma$ be a  subgroup of the group of translations of $R^n$ containing $Z^n$. Then the pseudogroup generated by the action of $\Gamma$ on $R^n$ is weakly isomorphic to the pseudogroup generated by the action of $\Gamma/Z^n$ on the torus $R^n/Z^n$. For instance the action of a subgroup $\Gamma$ of $R$ generated by two translations of amplitude $\alpha$ and $1$ is weakly isomorphic to the pseudogroup generated by a rotation of amplitude $\alpha$ modulo $1$  of the circle $R/Z$.

\noindent {\it Exercise}: Describe the group of weak automorphims of the pseudogroup $Q$ in the above example. Hint: it is isomorphic to the group of weak automorphisms of   $P$.

\noindent{\bf 1.6. The holonomy pseudogroup of a foliation}

Let $M$ be a differentiable manifold of dimension $n$ with a countable basis. Let $\Gamma_q$ be the etale groupoid of germs of local diffeomorphisms of $R^q$. A smooth foliation $\Cal F$ of codimension $q$ on $M$ can be given as follows. 
Choose an open covering $\Cal U =\{U_i \}_{i \in I}$ indexed by a  countable set $I$. For each $i \in I$ give a surjective smooth submersion $f_i: U_i 
\to V_i$ with connected fibers, where $V_i$ is an open subset of $R^q$. We can assume that the $V_i$ are disjoint.
We assume
 the following compatibility conditions: for each $(i,j) \in I \times I$ we give a continuous map $\gamma_{ij}:U_i \cap U_, \to \Gamma_q$ such that  $$f^x_i = \gamma_{ij}(x) f^x_j, \ \ \forall x \in U_i \cap U_j, $$
 where  $f^x_i$  and $f^x_j$ are the germs of $f_i$ and $f_j$ at $x$.  For $x \in U_i \cap U_j \cap U_k$, we have the cocycle relation :
 $$ \gamma_{ij}(x) \gamma_{jk}(x) = \gamma_{ik}(x).$$
 
 Note that $\gamma_{ii}(x) = f_i^x$ if the units of $\Gamma_q$ are identified with the points of $R^q$.

 Consider the topology on $M$ having as basis of open sets the intersections of the open sets of M with the fibers of the submersions $f_i$. The above compatibility conditions imply that the leaves of $\Cal F$ form a partition of $M$.   This topology will be called the leaf topology. A leaf of $\Cal F$ is a connected component of $M$ with respect to this finer topology. The inclusion of a leaf $L$ in $M$ is a smooth injective immersion.

Let $V$ be the disjoint union of the $V_i$. This is a smooth manifold of dimension $q$ with in general many connected components.

 By definition, the holonomy pseudogroup $H(\Cal F)$ (we should say related to the various choices above)  of the foliation $\Cal F$ is the pseudogroup of local diffeomorphisms of the manifold $V$ is  the pseudogroup whose associated  groupoid of germs is generated by the $g_{ij}(x)$.

It is easy to prove that the weak isomorphism class of the holonomy pseudogroups associated to two different choices is the same. So it is costumary to call the holonomy pseudogroup of $\Cal F$ any representative in its weak isomophism class.

As an example, consider the foliation of $R^2$ whose leaves are the fibers of the linear projection $p: R^2 \to R$ sending $(x,y)$ to $y$. Let $\Cal F$ be the foliation restricted
to the complement in $R^2$ of the origin $(0,0)$. The holonomy pseudogroup of this foliation is weakly isomorphic to the pseudogroup generated by the identity map of a non Hausdorff $1$-dimensional manifold. The projection $p$ induces an \'etale map from this manifold onto $R$. The fibers of this map consist of one point, except that the fiber above $0$ contains two points.

\input amstex
\documentstyle{amsppt}
\magnification=\magstep1

\noindent{\bf 1.7} {\bf Pullbacks of topological groupoids.} Let $\Cal G$ be a topological groupoid with space of objects $X$ and let $p:Y \to X$ be a continuous map of topological spaces. The pullback of $\Cal G$ by $p$ is the topological groupoid with space of objects $Y$, noted $p^{-1}(\Cal G)$, defined as follows. Its space of arrows  is the subspace  of the product $ Y \times \Cal G\times Y$ consisting of triples $(y_2, g,y_1)$ such that $p(y_1) =s(g)$ and $p(y_2)= t(g)$. The composition $({y_2}',g',{y_1}')(y_2,g,y_1)$ whenever defined is equal to $({y_2}',g'g,y_1)$. The map sending $(y_2,g,y_1)$ to $g$ is a homomorphism from $p^{-1}(\Cal G)$ to $\Cal G$.
 
 In the particular case where $Y$ is a subspace of $X$, the pullback of $\Cal G$ by the inclusion of $Y$ in $X$ is canonically isomorphic to the restriction of $\Cal G$ to $Y$.

 If $\Cal G$ is an \'etale groupoid, then $p^{-1}(\Cal G)$
is also an \'etale groupoid. This follows from the observation that the pullback of a homeomorphism $h: Z \to X$ by a continuous map $p: Y \to X$ is a homeomorphism $Y\times_X Z \to Y$.

Let $p:Y\to X$ be an \'etale map, $P$  a pseudogroup
of local homeomorphisms of $X$ and $\tilde P$ its associated groupoid of germs. Then $p^{-1}(\tilde P)$ is the groupoid of germs of the pseudogroup $p^{-1}(P)$ generated by the local homeomorphisms of $Y$ projecting by $p$ to elements of $P$.

\newpage

\heading {II. $(\Gamma \ltimes G)$-extensions of topological spaces}\endheading
Let $\Gamma$ be a dense subgroup of a simply connected Lie group $G$. Let $\Gamma \ltimes G$ be the topological groupoid given by the action of $\Gamma$ by left translations of $G$ (see I,1.2).  Roughly speaking a locally trivial $(\Gamma \ltimes G)$-extension of a topological space $X$ is a topological groupoid $\Cal G$ with a surjective map $p: \Cal G \to X$ such that each point $x$ has an open neighbourhood $U$ so that 
the subgroupoid $p^{-1}(U)$ is equivalent to the groupoid 
$\Gamma \ltimes( G \times U)$ associated to the action $\gamma.(g,u) \mapsto (\gamma g,u)$ of $\Gamma$ on $(G \times U)$.



\subheading{ 1. The groupoid $\Gamma \ltimes G$  and its group of automorphisms}

 Let $\Gamma$ be a dense subgroup of a connected Lie group $G$. We consider the groupoid $\Gamma \ltimes G$ defined by $\Gamma$ (with the discrete topology) acting by left translations on $G$. Note that it is the groupoid of germs associated to the pseudogroup of homeomorphisms of $G$ generated  by the left translations by the elements of $\Gamma$. We note $\rho: \Gamma \to G$ the inclusion of $\Gamma$ in $G$.

 Let $p_\Gamma$ (resp. $p_G$) be the projection of $\Gamma \ltimes G$ on the first factor (resp. the second one).

\noindent{\bf 1.1. Definition of $\alpha_\phi$ and $\phi_0$.}
For an automorphism $\phi$ of $\Gamma \ltimes G$, let $\alpha_\phi =p_\Gamma \circ \phi$. It is a map $(\gamma,g) \mapsto \alpha_\phi(\gamma,g)$ from $\Gamma \times G$ to $\Gamma$. As $G$ is assumed to be connected and $\Gamma$ has the discrete topology, this map does not depend on $g$ and will be noted $\gamma \mapsto \alpha_
\phi(\gamma)$

Let $\phi_0$ be the homeomorphism of $G$ defined by $\phi_0 = p_G \circ \phi$, 
so it is the homeomorphism of $G$  induced by $\phi$ on the space $G$ of objects of  $\Gamma\ltimes G$.
So $$\phi(\gamma,g) = (\alpha_\phi(\gamma),\phi_0(g)).$$
\proclaim{1.2. Lemma} For $\phi \in Aut(\Gamma \ltimes G)$, we have

\noindent 1) $\alpha_\phi$ is an automorphism of $\Gamma$.

\noindent 2) $\phi_0$ is $\alpha_\phi$-equivariant, i.e.
$$ \phi_0(\rho(\gamma)g) = \rho(\alpha_\phi(\gamma))\phi_0(g),\ \forall (\gamma,g) \in \Gamma \ltimes G.$$
3) For $\phi,\phi' \in Aut(\Gamma \ltimes G)$, then
$$\alpha_{\phi \circ \phi'}(\gamma) =\alpha_\phi(\alpha_{\phi'}(\gamma)),\ (\phi \circ \phi')_0(g)= \phi_0(\phi'_0(g)).$$
4) $\alpha_{\phi^{-1}}= (\alpha_\phi)^{-1}$ and
 ${(\phi}^{-1})_0 =(\phi_0)^{-1}$.
\endproclaim

\demo{Proof}
1) For two elements $(\gamma',g')$ and $(\gamma,g)$ of $\Gamma \ltimes G$, where $g' =\rho(\gamma)g$, their composition $(\gamma',g')(\gamma,g)$ is defined and is equal to $(\gamma'\gamma,g)$. So we have $\phi(\gamma',g')\phi(\gamma,g)= \phi(\gamma'\gamma,g)$. In the notations above this means that
 $$(\alpha_\phi(\gamma'),\phi_0(\rho(\gamma) g))(\alpha_\phi (\gamma),\phi_0(g)) =  (\alpha_\phi(\gamma')\alpha_\phi(\gamma),\phi_0(g))=(\alpha_\phi(\gamma'\gamma),\phi_0(g)).$$ Therefore $\alpha_\phi$ is a homomorphism from $\Gamma$ to $\Gamma$.
 
 We now prove that $\alpha_\phi$ is an automorphism of $\Gamma$. For $\phi',\phi \in Aut(\Gamma \ltimes G)$, we have $(\phi' \phi)(\gamma,g) = \phi'(\phi(\gamma,g)) =\phi'(\alpha_\phi(\gamma)),\phi_0(g)) =(\alpha_{\phi'}(\alpha_\phi(\gamma), \phi'_0(\phi_0(g))$. In particular for $\phi' =\phi^{-1}$, we see that $\alpha_{\phi^{-1}}$ is the inverse of $\alpha_\phi$.

2) For $\phi \in Aut(\Gamma \ltimes G)$ and $(\gamma,g) \in \Gamma \ltimes G$, we have  $\phi(\gamma,g) = (\alpha_\phi(\gamma), \phi_0(g))$. Then $\phi_0$ maps the target $\rho(\gamma)(g)$ of $(\gamma,g)$ to the target $\rho(\alpha_\phi)(\gamma)\phi_0(g)$ of $\phi(\gamma,g)$. This proves 2).

 3) and 4) We have $(\phi \circ \phi')(\gamma,g) =\phi(\alpha_{\phi'}(\gamma),\phi'_0(g)) =  (\alpha_\phi(\alpha_{\phi'}(\gamma),\phi_0(\phi'_0(g)) )$. This is equal to
$(\alpha_{\phi \circ \phi'}(\gamma), (\phi_0 \circ \phi'_0)(g)$. Comparing these two equalities gives the result. For 4), apply this result with $\phi' = \phi^{-1}$.

$\square$

\enddemo

\noindent{\bf 1.3. Definition of the group $Aut(\rho)$.} 

An element of $Aut(\rho)$ is an automorphism $\alpha$ of $\Gamma$ 
which is the restriction of a homeomorphism $\overline \alpha$ of $G$, more precisely:$$\overline \alpha(\rho(\gamma)) = \rho(\overline \alpha(\gamma)).$$
As we shall see, this homeomorphism $\overline \alpha$ is an
automorphism of $G$. Note that for an automophism $\phi\in Aut(\Gamma \ltimes G)$, then $\phi_0 \in Aut(\rho)$.

\proclaim {1.4. Proposition}
 1) For $\alpha \in Aut(\rho), \gamma \in \Gamma$ and $x\in G$ we have 
$$\overline \alpha (\rho(\gamma) x) = \rho(\alpha(\gamma)) \overline \alpha (x).$$
Moreover, $\overline \alpha$ is an automorphism of the Lie group $G$

2) The subgroup $Int(\rho)$ of $Aut(\rho)$  formed by the conjugate by the elements of $\Gamma$ is a normal subgroup of $Aut(\rho)$.

3) An element $\phi = (\alpha_\phi, \phi_0) \in Aut(\Gamma \ltimes G)$ is characterized by $\alpha_\phi$ and $\phi_0(1)$; more precisely $\phi_0(g) = \overline {\alpha_\phi}(g)\phi_0(1)$
\endproclaim
\demo{Proof} To simplify the notations we write $\gamma g$ instead of $\rho(\gamma)g$.

1) As $\Gamma$ is dense in $G$, there is a sequence $\gamma_n$  such that  $x= lim(\gamma_n)$. As $\overline \alpha$ is homeomorphism of $G$, we have 
$lim \ \overline\alpha(\gamma \gamma_n) = \alpha (\gamma)\  lim\ \overline\alpha(\gamma_n)= \alpha(\gamma)\overline \alpha(x)$.

2) For $\alpha \in Aut(\rho),\  \gamma,\gamma' \in \Gamma$, we have $(\alpha \circ ad(\gamma) \circ \alpha^{-1})(\gamma')= \alpha(\gamma \alpha^{-1}(\gamma')\gamma^{-1})=\alpha(\gamma)\gamma' \alpha(\gamma)^{-1} = ad(\alpha(\gamma))(\gamma')$.

3) Let $\phi = (\alpha_\phi,\phi_0) \in Aut(\Gamma \ltimes G)$. We have $\phi_0(\gamma x) =\alpha_\phi(\gamma)\phi_0(x)$ (see 2), Lemma 2.1). If $x=1$ this gives $\phi_0(\gamma)= \alpha_\phi(\gamma)\phi_0(1)$. Given $g\in G$, we can find a sequence $\gamma_n$ such that $g = lim \ \gamma_n$. Then $\phi_0(g) = lim \ \phi_0(\gamma_n) = \overline{\alpha_\phi}(g) \phi_0(1)$. $\square$

\enddemo

 Note that for $\phi\in Aut(\Gamma \ltimes G)$, the automorphism $\alpha_\phi$ belongs to  $Aut(\rho)$. 
 
 \proclaim{1.5. Corollary}(Definition of $Out(\rho)$). We have the exact sequence:
$$ 1 \to Int(\rho) \to Aut(\rho) \to Out(\rho) \to 1,$$
where $Out(\rho)$ is defined as the quotient of $Aut(\rho)$ by the normal subgroup $Int(\rho)$. We also have the exact sequence:
$$1 \to \Gamma_0 \to \Gamma  \to Aut(\rho) \to Out(\rho)\to 1,$$ where $\Gamma_0$ is the center of \  $\Gamma$ and $\Gamma \to Aut(\rho)$ maps $\gamma$ to $ad(\gamma)$.

In particular, if $G$ is abelian, then $Aut(\rho) = Out(\rho)$. 
 \endproclaim

The group $Aut(\rho)$ acts on this sequence by conjugation, more precisely an element $\alpha \in Aut(\rho)$ acts by conjugation on $Aut(\rho)$. It maps $\gamma \in \Gamma$ on $\overline \alpha (\gamma)$ and induces an action on $\Gamma_0$. It acts by the identity on $Out(\rho)$.

    \noindent{\bf 1.6. Definition of the group $Aut (\rho) \ltimes G.$} We have seen in proposition 1.4 that the group $Aut(\rho)$ acts by automorphisms of the Lie group $G$ as follows: the action of $\alpha \in Aut(\rho)$ on $g \in G$ is equal to $\overline \alpha(g)$.  The map $\alpha \mapsto \overline \alpha(g)$ from $Aut(\rho) $ to the group  of automorphisms of $G$ is a group homomorphism, i.e.
$$\overline{\alpha \alpha'}(g) =\overline\alpha (\overline{\alpha'}(g)),\ \forall g \in G, \ \forall \alpha, \alpha' \in Aut(\rho).$$ 
Therefore the composition in $Aut(\rho)\ltimes G$ is defined by
$$(\alpha,g)(\alpha',g') =(\alpha \alpha',g \overline \alpha(g') ).$$
Note that $(\alpha,g)^{-1}=(\alpha^{-1},\overline{\alpha}^{-1}(g^{-1}))$. 

Indeed $(\alpha,g)(\alpha^{-1},\overline{\alpha}^{-1}(g^{-1}))=(\alpha \alpha^{-1},\overline\alpha(g)\overline{\alpha}(g^{-1})) =(1,1)$.

The group $Aut(\rho) \ltimes G$ will be considered as a topological group, where $Aut(|\rho)$ is with the discrete topology, and $G$ with its usual topology.

\proclaim {1.7. Theorem} 1) The map $\psi : Aut(\rho) \ltimes G \to Aut(\Gamma \ltimes G)$ defined by 

 $$\psi(\alpha,g)(\gamma,x) = (\alpha (\gamma), \overline \alpha (x)g^{-1})$$
 
 is a group isomorphism.
 
 2) the map $\mu: \Gamma \to Aut(\rho) \ltimes G = Aut(\Gamma \ltimes G)$ defined by 
 $$\mu(\gamma) =(ad(\gamma),\rho(\gamma)^{-1})$$
 is a group homomorphism whose image is a normal subgroup. 
 \endproclaim

 \demo{Proof} 1). It follows from proposition 1.4, 3) that 
  the image of $\psi$ belongs to $Aut(\Gamma \ltimes G)$.
 Let us show that $\psi$ is a homomorphism, namely that
 $$\psi[(\alpha,g)(\alpha',g')] (\gamma,x)= \psi (\alpha,g)[\psi(\alpha',g')(\gamma,x)]
.$$

 The left hand side is equal to :
 
\noindent $\psi(\alpha \alpha',g \overline {\alpha} (g'))(\gamma,
x) = ((\alpha \alpha'(\gamma),\overline{\alpha \alpha'}(x) \overline \alpha(g')^{-1}g^{-1})$.    
  
    The right hand side is equal to :
    
 \noindent $\psi (\alpha,g)[\psi(\alpha',g')(\gamma,x)] =\psi(\alpha,g)[(\alpha'(\gamma),\overline \alpha'(x){g'}^{-1})] =(\alpha\alpha'(\gamma)\overline \alpha \overline \alpha'(x)  \overline {\alpha} ({g'}^{-1}) g^{-1}).$     
  
    This is equal to the left  hand side because $\overline{\alpha}({g'}^{-1}) =\overline{\alpha}(g')^{-1}$.

     The homomorphism $\psi$ is surjective; indeed given $\phi=(\alpha_\phi,\phi_0) \in Aut(\Gamma\ltimes G)$,  we have   $\psi(\alpha_\phi, \phi_0(1)^{-1})(\gamma,x) = (\alpha_\phi(\gamma),\overline{\alpha_\phi}(x)\phi_0(1)= (\alpha_\phi.\phi_0)$ by 	1.4. 3). It is also injective: if $(\alpha,g)$ belongs to the kernel of $\psi$, then we have $\psi(\alpha,g)(\gamma,x)= (\alpha(\gamma),xg^{-1}) =(\gamma,x)$; hence $\alpha =1$ and $g=1$.

2) The map $\mu$ is a group homomorphism, because for $\gamma,\gamma' \in\Gamma$, we have: 

\noindent  $(ad(\gamma \gamma'),\rho(\gamma 
  \gamma')^{-1})=(ad(\gamma),\rho(\gamma')^{-1})(ad(\gamma'),\rho(\gamma)^{-1})$; it is clear that $\mu$ is injective. Note that $\overline{ad(\gamma)}(g) =\gamma g \gamma^{-1}$.
  
 We now prove  that  $\mu(\Gamma)$ is a normal subgroup in $Aut(\Gamma \times G)$.
  We have: $$(\alpha,g) \mu(\gamma)(\alpha,g)^{-1} =(\alpha,g)(ad(\gamma),\rho(\gamma)^{-1})(\alpha^{-1}, \overline{\alpha^{-1}}(g^{-1}).$$
 This is equal to :
 $$[\alpha ad(\gamma),g\overline \alpha(\gamma^{-1})](\alpha^{-1},\overline{\alpha^{-1}}(g^{-1}))=$$

       $$(ad(\alpha(\gamma), g\overline\alpha(\gamma^{-1})\overline \alpha [(ad(\gamma)(\overline {\alpha^{-1}}(g^{-1}))] =$$
$$ = (ad(\alpha(\gamma), g\overline\alpha(\gamma^{-1})\overline \alpha (\gamma)g^{-1}\overline\alpha(\gamma^{-1})).
$$ This is equal to $ (ad(\alpha(\gamma)),\alpha(\gamma)^{-1})= \mu(\alpha(\gamma))$.
   $\square$

\enddemo

 \proclaim{1.8. Corollary } We have the  exact sequence of groups:
  $$1 \to \Gamma \to Aut(\Gamma  \ltimes G) \to Out(\Gamma \ltimes G) \to 1,$$
 where the first arrow is equal to $\mu$ and
$Out(\Gamma \ltimes G)$ is defined as the quotient of $Aut(\Gamma \ltimes G)$ by the image of $\mu$.

\endproclaim

\noindent{\bf 1.9. Remark.} Note that  $\mu(\Gamma)$ is not generally a closed subgroup of the Lie group $Aut(\Gamma \ltimes G)$; therefore, $Out(\Gamma \ltimes G)$ is not usually a Lie group. 

The subgroup $\mu(\Gamma)$ is closed in $Aut(\Gamma \ltimes G)$ in some particular cases, for instance when the center of $\Gamma$ is discrete (see the examples in part III). In that case 
$Out(\Gamma \ltimes G)$ with the quotient topology is a Lie group.

\proclaim{1.10. Proposition} The map $p: Aut(\Gamma \ltimes  G) = Aut(\rho) \ltimes G \to Aut(\rho)$ sending $(\alpha, g)$ to $\alpha$ can be considered as a fiber bundle with base space $Aut(\rho)$ and fiber $G$. The group $Aut(\rho)$ acts on this bundle as follows : the action of $\alpha'$ on $(\alpha,g) \in Aut(\rho) \ltimes G$ is equal to :
$\alpha'.(\alpha,g) = (\alpha'\alpha, \overline {\alpha'}(g)).$.

We have $\pi_i(Aut( \rho) \ltimes G, 1)$ is equal to $Aut(\rho)$ for $i=1$ and to $\pi_i(G,1)$ for $i >0$.

\endproclaim

\demo{Proof} We have to check that\ \  \ $(\alpha''\alpha').(\alpha, g) =\alpha''.(\alpha'.(\alpha,g))$.

The left side of this equality is equal to $(\alpha''\alpha'\alpha,  \overline {\alpha''\alpha'}(g))$. Its right side is equal to $\alpha''.( \alpha'\alpha,\overline {\alpha'}(g)= (\alpha''\alpha'\alpha, \overline {\alpha''\alpha'}(g))$. 

To check the last sentence, use the homotopy exact sequence of this bundle.
$\square$
\enddemo
\noindent {\bf 1.11. Some examples}
 
 1) Assume that $\Gamma$ is the group $G$ with the discrete topology and that $\rho$ is the natural inclusion. Then $Aut(\rho)$ is the group of automorphisms of $G$ (as a Lie group) with the discrete topology.  
 
 2) Let $\rho(\Gamma)$ be the subgroup $Z[\sqrt 2]$ of R. An element of $Aut(\rho)$ is an automorphism $\lambda$ of $\rho(\Gamma)$ given by the multiplication by an element of the form
 $a+ b\sqrt 2$, where $a,b \in Z$. Therefore $\lambda$ can be represented  with repect to the basis  $(1,\sqrt 2)$ by a matrix of the   form:
 
 $$M =\left( \matrix
 a & 2b\\
 b & a
 \endmatrix \right )
 $$
 As $\lambda$ is invertible, we have $a^2 - 2b^2 = \pm 1$.
  Therefore $Aut(\rho) \in Sl_2(Z)$.
  
  The elements  $a+ b\sqrt 2$, with $a,b \in Z$ and $a^2-2b^2 =\pm 1$, form a subgroup generated by $-1$ and $1+\sqrt 2$ of the group of units of $Q(\sqrt 2)$, we have $Aut(\rho) = Z_2 + Z$. 
   
    Note that for each dense subgroup $\rho(\Gamma)$ of $R$, then 
  $Aut(\rho)$ contains the cylic group of order 2 represented by the symmetry of $R$ fixing $0$.
  
  For more examples, see Part III,2.
  
3) Let $\rho: \Gamma \to G$ be the inclusion of $\Gamma$ as a dense subgroup of $G$. Then $Aut(\rho \times \rho) =Z_2
 \ltimes (Aut(\rho) \times Aut(\rho))$, where $Z_2$ acts by exchange of the factors. Similarly $Out(\rho \times \rho) = Z_2 \ltimes(Out(\rho) \times Out(\rho))$.
 
 4) Consider the case where $G =SU_2 =S^3$ and $\rho:\Gamma \to SU_2$. The group of autmorphisms of 
$G$ is equal to $Z_2 \times G$, where the elements of $G$ acts by conjugation (note that the center of $G$ is trivial) and $Z_2$ by complex conjugation.  This is also true for $SU(n)$ where $n>1$. Hence $Aut(\rho) =Z_2 \times \Gamma$.
 .

 5) It is an interesting exercise to generalize the results of this
  paragraph 1 to the case where the inclusion of $\Gamma$ in $G$ is replaced by a group homomorphism $\rho: \Gamma \to G$, where $\Gamma$ is a group with the discrete topology  and the image of $\rho$ is a dense subgroup of $G$.

\newpage
 
\noindent{\bf 2. $(\Gamma \ltimes G)$-extensions of topological spaces.}
 
   Throughout  this section we fix a dense subgroup $\Gamma$ of $G$; the inclusion of $\Gamma$ in $G$ is noted $\rho$. 
  We define the notion of a  $(\Gamma\ltimes G)$-extension of a topological space $X$, of the set of isomorphism classes of such extensions and the set of equivalence classes of such extensions. 
  
  The determination of the set of homotopy classes
of equivalence classes of $(\Gamma \ltimes G)$-extensions of $X$  is the main purpose of this paper.

 \noindent{\bf 2.1. Definition of  $(\Gamma \ltimes G)$-extensions.}
 
 For a topological space $U$, the groupoid $(\Gamma \ltimes G) \times U$ is the groupoid associated to the action of $\Gamma$ on $G\times U$, where  $\gamma\in\Gamma$  maps $(g,u)$ to $(\rho(\gamma) g,u)$. This groupoid will be considered as a trivial $(\Gamma \ltimes G)$-extension of $U$.

 A $(\Gamma \ltimes G)$-extension of a topological space $X$ is a topological groupoid $\Cal G$ of  germs (see I,1.4)  with a surjective continuous map $p:\Cal G \to X$. We assume that each point $x \in X$ has an open neighbourhood $U$ such that there is a weak isomorphism (see I,1.5) noted $E_U: p^{-1}(U) \to (\Gamma \ltimes G) \times U$  commuting with the natural projections to $U$. Moreover the following compatibility condition must be satisfied : if $U_1$ and $U_2$ are such  open sets of $X$ such that $U_1\cap U_2$ is non empty and if $E_1: p^{-1}(U_1) \to (\Gamma \ltimes G) \times U_1$  and $E_2: p^{-1}(U_2) \to (\Gamma \ltimes G) \times U_2$ are weak isomorphisms, then the  restriction to $U_1\cap U_2$ of  $E_1\circ {E_2}^{-1}$ is an isomorphism of $(\Gamma \ltimes G) \times  (U_1\cap U_2)$ of the form  $(\gamma,g,u) \mapsto (\phi_{12}(\gamma,g),u)$ where $\phi_{12}$ is a continuous map $U_1\cap U_2 \to Aut(\Gamma \ltimes G)$.

 A weak isomorpism of $(\Gamma \ltimes G)$-extensions from $p:\Cal G \to X$ to $p':\Cal G' \to X$   is a weak isomorphism from $\Cal G$ to $\Cal G'$ commuting with $p$ and $p'$. More explicitely, each point of $X$ is contained in an open set $U$ such that there are weak isomorphims $E_U: p^{-1}(U) \to (\Gamma \ltimes G) \times U$ and $E'_U:{p'}^{-1}(U) \to (\Gamma \ltimes G) \times U$. Then $E'_U \circ E^{-1}_U$ should  be of the form $(\gamma, g,u) \mapsto (\phi_u(\gamma,g),u)$, where $\phi_u$ is a continuous map from $U$ to $Aut(\Gamma \ltimes G)$.

\noindent {\bf 2.2. Definition of a
  $(\Gamma\ltimes G)$-extension  
 over an open cover}
 
 Let $\Cal U = \{U_i\}_{i \in I} $ be an open cover of $X$. Consider the disjoint union of the $U_i$, i.e, the space  $\bigcup_{i \in I}(U_i,i)$.   A $(\Gamma \ltimes G)$-extension over $\Cal U$ is an extension $p: \Cal G \to  \bigcup_{i\in I}(U_i,i)$ satisfying the following conditions.
 
 Its space of objects is the union of the $G\times (U_i,i)$. The restriction $\Cal G_i$ of $\Cal G$ to $G\times (U_i,i)$   is the trivial $(\Gamma \ltimes G)$-extension by $(U_i,i)$, so $\Cal G_i = (\Gamma \ltimes G) \times (U_i,i)$.
 
 For $(i,j)\in I\times I$, the restriction $\Cal G_{ij}$ of $\Cal G$ to  $(\Gamma \ltimes G) \times (U_i \cap U_j,j) \in \Cal G_j$ is of the form $(\gamma,g,(u.j)) \mapsto (\phi_{ij}(\gamma,g)(u,i))$, where $\phi_{ij}$ is a continuous map from $U_i \cap U_j$ to $Aut(\Gamma \ltimes G)$.

We assume that $\phi_{ii}$ is the identity and the cocycle relation:
 $$\phi_{ik}(u) =\phi_{ij}(u) \phi_{jk}(u),\ \ \forall u \in U_i \cap U_j \cap U_k.$$    It follows  that  $\phi_{ji} = {\phi_{ij}}^{-1}$. 
 
 \noindent {\bf 2.3. Model for a $(\Gamma\ltimes G)$-extension.} 
 
 Let $p:\Cal G \to X$ be a $(\Gamma \ltimes G)$-extension of $X$    We can find a fine enough open cover $\Cal U = \bigcup_{i \in I} U_i$ such that, for each $i \in U_i$, there is a weak isomorphism  $E_i: p^{-1}(U_i) \to (\Gamma \ltimes G) \times U_i$ commuting with the projection to $U_i$. The compatibility condition mentionned in 2.1 implies that for $u \in U_i \cap U_j$ then $E_j \circ {E_i}^{-1}$ is given by a continuous map $\phi_{ji}: U_j \cap U_i \to Aut(\Gamma \ltimes G)$. 
 Therefore we see that we get a weak isomorphism from $p:\Cal G \to X$ to a $(\Gamma \ltimes G)$- extension  over $\Cal U$.

 If we had chosen instead of $E_i$ another weak isomorphism $E'_i :p^{-1}(U_i) \to (\Gamma \ltimes G) \times U_i $, then $E'_i \circ E_i^{-1}$ would be an  isomorphism of $(\Gamma\ltimes G) \times  U_i$ of the form $((\gamma,g),u) \mapsto (\phi_i(u)(\gamma,g),u)$, where $\phi_i:U_i \to Aut(\Gamma \ltimes G)$. The corresponding $\phi'_{ij}(u)$ would be equal to $\phi_i
 (u)\phi_{ij}(u) \phi_j^{-1}(u)$ for $u \in U_i \cap U_j$.

\noindent{\bf 2.4. Definition of
$H^0(X,\Cal A)$ and $H^1(X,\Cal A)$.}

In this section we recall  basic notions about the cohomology of sheaves of groups which are not necessary commutative. The reference is Serre  [9] Chapitre I, pp. 199-222.

Let $\Cal A$ be a sheaf of groups over $X$. The stalk of $\Cal A$ over $x \in X$ is a group which is not abelian in general. For instance, $\Cal A$ could be the sheaf of germs of continuous maps from open sets of $X$ to a topological group.

Let $\Cal U = \{U_i\}_{i \in I}$ be an open cover of $X$ indexed by a set $I$. Let $\Cal C ^0(\Cal U, \Cal A)$ be the set whose elements associate to each $i\in I$ a continuous section $\phi_i$ of $\Cal A$ over $U_i$. Similarly let $\Cal C^1(\Cal U, \Cal A)$ be the set associating to each couple $(i,j) \in I \times I$ such that $U_i \cap U_j$ is non empty a continuous section $\phi_{ij}: U_i \cap U_j \to \Cal A$ .

We first define $H^0(\Cal U,\Cal A)$. It is the subset of $\Cal C^0(\Cal U,\Cal A)$ made up of cochains $\phi$ such that $\phi_i(u) =\phi_j(u)$ for each $u \in U_i \cap U_j$. So the maps $u\mapsto \phi_i (u)$ for each $i\in I$ fit together to define a continuous section of $\Cal A$. Therefore $H^0(\Cal U,\Cal A)$ can be identified to the group of continuous sections of $\Cal A$; it is independent of the open cover $\Cal U$ and will be noted $H^0(X,\Cal A)$.

We now define $H^1(\Cal U,\Cal A)$.  Let $\Cal Z^1(\Cal U, \Cal A)$ be the subset of $\Cal C^1(\Cal U, \Cal A)$ formed by the $1$-cocycles, i.e. the elements such that $$\phi_{ij}(u)\phi_{jk}(u) =\phi_{ik}(u),\ \ \forall u \in U_i \cap U_j \cap U_k, \ \ (i,j,k) \in I^3.$$
Note that this implies that $\phi_{ii}(u)=1_u$ and $\phi_{ij}(u) = {\phi_{ji}}^{-1}(u)$.

 For $\phi \in \Cal C^0(\Cal U,\Cal A)$,
 let $\delta :\Cal C^0 \to \Cal C^1$ be the map defined by :
  $$(\delta\phi)_{ij}(u) =\phi_i(u)\phi_{ij}(u)\phi_j^{-1}(u)=\phi'_{ij}, \  \ \forall u \in U_i\cap U_j.$$

 We say that $\phi'_{ij} =\phi_i\phi_{ij}\phi_j^{-1}$ is related to $\phi_{ij}$ by the elements $\phi_i ,\phi_j\in \Cal C^0(\Cal U,\Cal A)$. This defines an equivalence relation in $\Cal Z^1(\Cal U,\Cal A)$; indeed, $\phi_{ij}$ is related to $\phi'_{ij}$ by $\phi_i^{-1}$  $\phi_j^{-1}$  and if  $\phi''_{ij}$ is related to $\phi'_{ij}$ by $\phi'_i$, $\phi'_j$, then $\phi''_{ij}$ is related to $\phi_{ij}$ by $\phi'_i\phi_i,\phi'_j\phi_j$. Also $\phi_{ij}$ is related to itself by the $0$-cochain associating to each $u \in U_i$ the unit element $1_u$  and to each $u\in U_j$ the unit $1_u$.
 The set of equivalence classes in $\Cal Z^1(\Cal U,\Cal A$ is noted  $H^1(\Cal U,\Cal A)$.
 
 Let $\Cal V =\{V_l\}_{l \in L}$ be an open cover of $X$ indexed 
by a set $L$, which is finer than $\Cal U$. This means that there is a map $\tau:L \to I$ such that $V_l \subset U_{\tau(l)}$ for all $l \in L$. The map $\tau$ induces, for $r=0,1$, maps $\tau^r:H^r(\Cal U,\Cal A) \to H^r(\Cal V, \Cal A)$ defined as follows.

For $r=0$, an element of $H^0(\Cal U,\Cal A)$ is represented by a $0$-
cocycle $\{\phi_i\}$, i.e such that $\phi_i(u)= \phi_j(u)$ for all $u \in U_i \cap U_j$. Then its image by $\tau^0$ is the element of $H^0(\Cal V,\Cal A)$ represented by the $0$-cocycle $\{\phi_l\}$ defined  for $v \in V_l$ by $\phi_l(v) =\phi_{\tau(l)}(v)$. This map defines a group homomorphism which does not depend on the choice of $\tau$. Moreover it is an isomorphism because both groups are naturaly isomorphic to the group $H^0(X,\Cal A)$ of global sections of $\Cal A$.

Similarly, for $r=1$, an element of $H^1(\Cal U, \Cal A)$ represented by the cocycle $\{\phi_{ij}\}$ is mapped by $\tau^1$ to the element represended by the cocycle $\{\phi_{lm}\}$ defined for $v \in V_l \cap V_m$ by $\phi_{lm}(v) =  \phi_{\tau(l)\tau(m)}$. As above, $\tau^1: H^1(\Cal U,\Cal A) \to H^1(\Cal V,\Cal A)$ does not depend on the choice of $\tau$; it is a map of sets preserving the distinguished base points. Moreover, it is injective (cf. [9]).

We  consider only open covers finer than $\Cal U$ 
(more precisely equivalence classes of such covers in the sense of Serre ([9], p.215). They form an ordered filtered set
and  the set $H^1(X, \Cal A)$ is defined as the inductive limit of the sets $H^1(\Cal U,\Cal A)$ (see [9]).

\proclaim{2.5. Theorem}Let $\Cal A$ be the sheaf $\tilde{Aut}(\Gamma \ltimes G))$ over  $X$ of germs of continuous maps of open sets of $X$  to the topological group $Aut(\Gamma \ltimes G)$.

 The set of weak isomorphism classes of $(\Gamma \ltimes G)$-extensions $p: \Cal G \to X$ of a topological space $X$ is in bijection with the set $H^1(X, \tilde{Aut}(\Gamma\ltimes G))$.
\endproclaim

 \demo  {Proof} Let $p: \Cal G \to X$ be a $(\Gamma \ltimes G)$-extension as above.
Let us choose an open cover $\Cal U = \{U_i\}_{i \in I}$ of $X$ indexed by a set $I$ such that for each $i \in I$ there is a weak isomorphism $E_i: p^{-1}(U_i) \to (\Gamma \ltimes G )\times U_i$. Given $j \in I$ such that $U_i \cap U_j \neq \emptyset$, then the restriction of $E_i \circ E_j^{-1}$ to $(\Gamma \ltimes G )\times (U_i \cap U_j)$ is of the form $(\gamma,g,u) \mapsto (\phi_{ij}(u)(\gamma,g)$, where $\phi_{ij} : U_i \cap U_j \to Aut(\Gamma\ltimes G)$ is a continuous map.  For $i,j,k \in I$ such that $U_i \cap U_j \cap U_k \neq \emptyset$, we have then the cocycle relation
$$  \phi_{ij}(u)\phi_{jk}(u) = \phi_{ik}(u), \ \forall u \in U_i \cap U_j \cap U_k. \tag  1$$

 The relation (1) defines a 1-cocycle over $\Cal U$ with value in the sheaf $\Cal A$, so an element of the set $Z^1(\Cal U,\Cal A)$. 
Let $\phi \in H^1(X, \tilde {Aut}(\Gamma \ltimes G))$ be the cohomology class represented by the $1$-cocycle $\phi_{ij}$.

Consider $p: \Cal G \to X$ and $p': \Cal G' \to X$  two $(\Gamma \ltimes G)$-extensions such that there is a weak isomorphism 
$w : \Cal G \to \Cal G'$ commuting with the projections to $X$. Choose a fine enough open cover $(U)_{i\in I}$ of $X$ such that there are for each $i \in I$ weak isomorphisms $E_i: p^{-1}(U_i) \to (\Gamma \ltimes G) \times U_i$ and $E'_i:{p'}^{-1}(U_i) \to (\Gamma \ltimes G) \times U_i$.

 Then $E'_i\circ E_i^{-1}$ is of the form $(\gamma,g,u) \mapsto (\phi_i (u)(\gamma,g),u)$ and $\{\phi_i\}_{i \in I}$ can be considered as an element of the group $C^0(\Cal U, \Cal A)$ of $0$-cochains over $\Cal U$ with value in $\Cal A$. 
We have $\phi'_{ij}(u) = \phi_i (u) \phi_{ij}(u)  \phi_j^{-1}(u)$ for $u\in U_i \cap U_j$. So the cocycles $\phi_{ij}$ and $\phi'_{ij}$ define the same element of $H^1(\Cal U,\Cal A)$.

 Therefore we have associated to the extension $p:\Cal G \to X$ a well defined element  $ \phi \in H^1(X,\Cal A)= H^1(X, \tilde{Aut}(\Gamma \ltimes G))$.

Conversely, given a $1$-cocycle $\{\phi_{ij}\} \in Z^1(\Cal U,\Cal A)$ representing a given element of $H^1(X,\Cal A)$, we can reconstruct a $(\Gamma \ltimes G)$-extension $p':\Cal G' \to X$ as in 2.3.
 $\square$
\enddemo
\noindent {\bf 2.6. Definition of equivalences of $(\Gamma \ltimes G)$-extensions.}

We consider now the set of equivalence classes of $(\Gamma \ltimes G)$-extensions of $X$. In 1.7 we have defined an injective homomorphism $\mu:\Gamma \to Aut(\Gamma \ltimes G)$. The image $\mu(\Gamma)$ is a normal subgroup of $Aut(\Gamma \ltimes G)$, but this subgroup is not closed in general, as noticed before. Therefore the quotient $Out(\Gamma \ltimes G) =Aut(\Gamma \ltimes G)/\mu(\Gamma)$ is not a Lie group in general.

 An equivalence of $(\Gamma \ltimes G)$-extensions from $p:\Cal G \to X$ to $p':\Cal G' \to X$   is given by a weak isomorphism $\epsilon: \Cal G \to \Cal G'$ commuting with the projections, such that, each point of $X$ is contained in an open set $U$ and there are weak isomorphims $E_U: p^{-1}(U) \to (\Gamma \ltimes G) \times U$ and $E'_U:{p'}^{-1}(U) \to (\Gamma \ltimes G) \times U$ in a way that  $E'_U \circ E^{-1}_U$ are of the form $(\gamma, g,u) \mapsto (\phi_u(\gamma,g),u)$, where $\phi_u$ is a continuous map from $U$ to $Out(\Gamma \ltimes G)$.

Let $\tilde {\Gamma}$ be the sheaf over $X$  of  germs of continuous maps from open sets of $X$ to the discrete group $\Gamma$. The homomorphism $\mu$ induces an injective map $\tilde {\mu}:\tilde {\Gamma} \to  \tilde {Aut}(\Gamma \ltimes G)$ of sheaves over $X$. The quotient of $\tilde {Aut}(\Gamma \ltimes G)$ by the subsheaf  $\tilde{\mu}(\tilde {\Gamma})$ is noted $\tilde {Out}(\Gamma \ltimes G)$. We have the exact sequence of sheaves over $X$:
$$1 \to \tilde \Gamma \to \tilde {Aut}(\Gamma \ltimes G) \to \tilde{Out}(\Gamma \ltimes G) \to 1,$$
and the projection $ \tilde {Aut}(\Gamma \ltimes G) \to \tilde {Out} (\Gamma \ltimes G)
$ is noted $\tilde \psi$.

\proclaim {2.7. Theorem} The set of equivalence classes of $(\Gamma \ltimes G)$- extensions of $X$ is equal to $H^1(X, \tilde {Out}(\Gamma \times G))$. \endproclaim

\demo{Proof} With the notations of the proof of theorem 2.5,
let $p:\Cal G \to X$ be a $(\Gamma \ltimes G)$-extension characterized by a $1$-cocycle $\phi_{ij} \in Z^1(\Cal U,\tilde {Aut} (\Gamma \ltimes G))$. Let $\psi_{ij} \in Z^1(\Cal U, \tilde{Out}(\Gamma \ltimes G))$ be its image induced by $\tilde \psi$.
Let $\psi \in H^1(X,\tilde{Out}(\Gamma \ltimes G))$ be the cohomology class represented by the cocycle $\psi_{ij}$.
As in the  proof of theorem 2.5, we see that the cohomology class $\psi$ does  not depend on the various choices in its construction and characterizes the equivalence class of $p: \Cal G \to X$.

Consider now another $(\Gamma \ltimes G)$-extension $p':\Cal G' \to X$ and an equivalence $\epsilon: \Cal G \to \Cal G'$ commuting with the projections to $X$. With the notations of the proof of theorem 2.5, we can find an open cover $(U)_{i \in I}$ of $X$ and 
 weak isomorphisms $E_i:p^{-1} (U_i) \to (\Gamma \ltimes G) \times U_i$ and $E'_j: {p'}^{-1}(U_j) \to (\Gamma \ltimes G) \times U_j$. The composition of $E'_i \circ {E_j}^{-1}$ with $\psi$ is noted $\psi_{ij}$.
It is a map from $U_i \cap U_j$  to $\tilde {Out} (\Gamma \ltimes G)$ which is a $1$-cocycle $Z^1(\Cal U, \tilde {Out} (\Gamma \ltimes G))$ representing an element of $H^1(X, \tilde {Out}(\Gamma \ltimes G))$ equal to $\psi$. $\square$.

\enddemo

\noindent{\bf 2.8. Functorial properties.}

Let $p: \Cal G \to X$ be a $(\Gamma \ltimes G)$-extension of $X$. Its pullback by a continuous map $f: X' \to X$ is the $(\Gamma \ltimes G)$-extension $p': \Cal G' \to X'$ of $X'$ defined as follows. As a topological space $\Cal G'$ is the subspace of $ \Cal G \times X'$ formed by the pairs $(g,x')$ such that $f(x') = p(g)$. The composition is defined by $(g_1,x')(g_2,x') =(g_1g_2,x')$ and $p'(g,x') = x'$. The groupoid $\Cal G'$ will be also denoted $f^{-1} \Cal G$. If $f': X''\to X'$ is a continuous map of topological spaces, then $(f \circ f')^{-1} \Cal G = {f'}^{-1}(f^{-1}\Cal G)$.
Also, if $\overline {\Cal G}$ is equivalent to $\Cal G$, then $f^{-1} \overline {\Cal G}$ is equivalent to $f^{-1}\Cal G$. If $f$ is the inclusion  of a subspace $X'$ of $X$, then $f^{-1}(\Cal G)$ is called the restriction  of $\Cal G$ to $X'$ and can be noted $\Cal G|_{X'}$

The map $\Cal G \mapsto f^{-1} \Cal G$ can be considered as a contravariant functor in the following sense:
if $E : \Cal G_1 \to \Cal G_2$ is a weak  isomorphism of $(\Gamma \ltimes G)$-extensions over $X$, then $f^{-1}E : = f^{-1}\Cal G_2 \to f^{-1}\Cal G_1$ is a weak isomorphism  over $X$ mapping $(g_2,f(x))$ to $E(g_1,x)$ is an equivalence.

On the other hand,  a continuous map $f:X' \to X$ induces a map $H^1(X, \Cal A) \to H^1(X', f^{-1}(\Cal A))$ called $H^1(f^{-1})$, and $f^{-1}(\Cal A)$ is the sheaf of germs of continuous maps from open sets of $X'$ to $Aut(\Gamma \ltimes G)$. According to theorem 2.5, an element of $H^1(X,\Cal A)$ represents   a weak isomorphism class of an extension $\Cal G$ and its image by $H^1(f^{-1})$ represents a weak isomorphism class of $f^{-1}\Cal G$.

Similarly, for equivalence classes of $(\Gamma \ltimes G)$-extensions, one has to replace the sheaf $\Cal A=  \tilde {Aut}(\Gamma \ltimes G)$ by the sheaf $ \tilde {Out}(\Gamma \ltimes G)$.

\noindent{\bf 2.9. Homotopy of $(\Gamma \ltimes G)$-extensions.}

Let $p_0 :\Cal G_0 \to X$ and $p_1:\Cal G_1 \to X$ be two $(\Gamma \ltimes G)$-extensions of $X$.
A homotopy from $p_0$ to $p_1$ is a $(\Gamma \ltimes G)$-extension 
$p: \Cal G \to X\times [0,1]$ and weak isomorphisms $E_i:\Cal G|_{X \times i} \to \Cal G_i$, where $i = 0$ or $1$. The relation of homotopy between extensions of $X$ is an equivalence relation. The pullback of an homotopy is an homotopy. We note $H(X)$ the set of homotopy classes of $(\Gamma \ltimes G)$-extensions of $X$. We get a contravariant functor noted $H$ from the category of topological spaces to the category of $(\Gamma \ltimes G)$-extensions by associating to a continuous map $f:X' \to X$ the pullback $f^{-1}: H(X) \to H(X')$. The functor $H$ will be also noted $H_{Aut (\Gamma  \ltimes G)}$. This notation will be justified below.

Similarly, an homotopy from the equivalence class of $p_0:\Cal G \to X$ to the equivalence class of $p_1: \Cal G_2 \to X$ is an equivalence class of $(\Gamma \ltimes G)$-extension of $X \times [0,1]$ together with equivalences $E_i$  from the restrictions of this class to $X \times i$ to $p_i:\Cal G \to X$, for $i= 0$ or $1$.

\noindent{\bf 2.10. Homotopy of $(\Gamma \ltimes G)$-extensions preserving a base point.}
Let $*$ be a base point in $X$. For a sheaf $\Cal A$ over  $X$, let $\Cal A_*$ be  the subsheaf of $\Cal A$ whose stalk above $*$ consists only of the neutral element.

 . 

With the notations of 2.9,
assume that $\Cal G_0|* = \Cal G_1|*$. A homotopy from $p_0$ to   $p_1$ preserving the base point $*$ is a $(\Gamma \ltimes G)$-extension $p : \Cal G \to X \times [0,1]$ as above such that the restriction of $\Cal G$
to $* \times [0,1]$ is the trivial extension of $\Cal G_0|* =\Cal G_1|*$ by $[0,1]$. We note $H(X,*)$ the set of homotopy of weak equivalemce classes of $(\Gamma\ltimes G)$- extensions of $X$   with base point $*$. We can consider $H(X,* )$ as a contravariant functor from the category whose objects are base point preserving continuous maps of pointed topological spaces to the category whose objects are the homotopy classes of weak equivalence classes preserving a base point.

The analogue of theorem 2.5 and 2.7 is the following:

\proclaim{Theorem 2.11}1) The set of homotopy classes of weak isomorphisms classes of $(\Gamma    \ltimes G)$-extensions of $X$ preserving  the base point $*$ is in bijection with the set $H^1(X, \tilde{Aut}(\Gamma \ltimes G)_*)$.

2) The set of homotopy classes of equivalence classes of $(\Gamma \ltimes G)$-extensions of $X$ preserving the base point $*$ is in bijection with the set $H^1(X, \tilde{Out}(\Gamma \ltimes G)_*)$.
\endproclaim 

\noindent{\bf 2.12. The smooth case.} Let $X$ be a paracompact smooth manifold. We use the word smooth to mean $C^\infty$-differentiable. A smooth $(\Gamma \ltimes G)$-extension  of $X$ is a $(\Gamma \ltimes G)$-extension $p: \Cal G \to X$ where $\Cal G$ is a groupoid of germs associated to a differentiable groupoid (cf. I.5.). In the definition 2.1 of a $(\Gamma \ltimes G)$-extension, we replace everywhere continuous by smooth.

Theorems 2.5 and 2.6 generalize to the case of smooth extensions with the same proof if we replace the sheaves of germs of continuous maps by the sheaves of germs of differentable maps.

Similarly if $X$ is a real analytic manifold, we have to replace everywhere continuous by real analytic. Also if $G$ is a complex Lie group (for instance $SU(2)$), we get the notion of complex analytic $(\Gamma \ltimes G)$-extensions by replacing everywhere continuous by holomorphic. The analogue of theorems 2.5 and 2.6 are also valid in those cases.

One has to precise the notion of homotopy in those cases. In the notations of 2.8, we assume that the weak isomorphisms $E_i :\Cal G_{|X \times i}\to \Cal G_i$ for $i=0$ or $1$ are the identity in a nbhd of $X\times i$. With this condition, the relation of homotopy is an equivalence relation. 

\proclaim{Theorem 2.12} Two smooth $(\Gamma \ltimes G)$-extensions which are topologically homotopic are smoothly homotopic. \endproclaim
For the proof, see [12].

\newpage
\noindent {\bf 3. Classifying spaces. }

In this section, we recall the main features of classifying spaces, namely their definition, their main properties and their construction. Then, we shall apply those generalities to the specific examples related to the purpose of our paper.

\noindent {\bf 3.1. The classifying space of a topological groupoid $G$.}

{\bf 3.1.1. $G$-principal bundles.}
In [3] , the notion of $G$-principal bundle over a topological space $X$ is defined as a natural generalization of the classical case where $G$ is a topological group (in loc.cit., $G$ is denoted by $\Gamma$ and a $G$-principal bundle
over $X$ is called a $\Gamma$-structure over $X$).

A $G$-principal bundle  $p: E \to X$ with total space $E$,  base space a topological space $X$ and projection $p$ is a locally trivial bundle with fiber $G$ endowed with a right action of $G$
acting simply transitively on each fiber.

There is a natural bijection between the set of isomorphism classes of $G$-principal bundles over $X$ and the set
$H^1(X, \tilde G)$, where $�\tilde G$ is the sheaf over $X$ of germs of continuous maps from open sets of $X$ to $G$. Indeed, consider a $G$-principal bundle $E$ over $X$ and let $\{U_i \}_{i \in I}$ be an open cover of $X$ such that $p^{-1}(U_i)$ has a continuous section $\sigma_i$. For $x \in U_i \cap U_j$, there is a unique continuous map $g_{ij}(x)$ such that 
$\sigma_j(x) = \sigma_i (x)g_{ij}(x)$. Then, $g_{ij}$ is a $1$-cocycle representing an elements of $H^1(X,\tilde G)$. Conversely, given such a cohomology class, one can reconstruct the $G$-principal bundle $E$.

Two $G$-principal bundles $p_0: E_0 \to X$ and $p_1:E_1 \to X$ are homotopic if there is a $G$-principal bundle $E
\to X \times [0,1]$ such that $ E$ restricted to $(X , 0)$ and $(X ,1)$ is equal to $E_0$ and $E_1$. Note that, if $G$ is a topological group, then, two  $G$-principal bundles which are homotopic, are isomorphic. This is not true in general for a topological groupoid.

Let $f: X' \to X$ be a continuous map of topological spaces. Then $f^*E$ is the $G$-principal bundle over $X'$ defined by the cocycle $g'_{ij} (x') = g_{ij}(f(x'))$ over the open cover $f^{-1}\{U_i\}_{i\in I}$ of $X'$.

Let $H_ G $ be the functor associating to each paracompact space $X$ the set $H_ G (X)$ of homotopy classes of principal $ G$-bundles over $X$. A continuous map $f: X' \to X$ of paracompact spaces induces a map $H_G(f) : H_G(X) \to H_G(X')$ associating to a principal $G$-bundle over $X$ its pulback by $f$. Also if $h: G_1 \to G_2$ is a continuous homomorphism of topological groupoids, it induces for each $X$ a map $H(h):H_{G_1}(X) \to H_{G_2}(X)$. So, $H_G(X)$ is a functor contravariant in $X$ and covariant in $G$.

{\bf 3.1.2. The classifying space of a topological groupoid.}
The joint construction of Milnor for the classifying space of a topological group has been generalized to the case of a topological groupoid $G$ by Lor and Buffet in 1970. It provides a canonical model for a universal principal $G$-bundle $p:EG \to BG$.

Let  $Y$ be the space of units of the topological groupoid $ G$.

Let us  recall this construction of $BG$ (see [3] with other notations). The elements of $G$ will be noted $g$ and the units (or the objects of $G$) will be noted $y$. The space $BG$ is the quotient of the space $EG$ considered as a principal $G$-bundle with base space $BG$, quotient of $EG$ by a right action of $G$.
 
 The elements of $EG$ are infinite sequences $(t_0,g_0,t_1,g_1,...)$, where all the $t_i$ are real numbers in $[0,1]$ which are all zero except for a finite number of indices and $\sum t_i = 1$. The elements $g_i$ belong to $G$ and have the same right unit (or source). Two such sequences $(t_0,g_0,t_1,g_1,...)$ and $(t'_0,g'_0,t_1,g_1,...)$ are equivalent if $t_i = t'_i$ for each $i$ and $g_i =g'_i$ if $t_i \neq 0$; hence, if $t_i = 0$ , $g'_i$  might be any element of $G$ with the same target
as $g_i$.  
This justifies the notation $(t_0g_0, t_1g_1,...)$ for the class of $(t_0,g_0,t_1,g_1,...)$.

 $G$ acts naturally on the right of $EG$ as follows. Given $g \in G$ with target $y$, its action on
$(t_0g_0,t_1g_1,...)$ in $EG$, where the  source of each $g_i$ is $y$, is equal to  $(t_0g_0g,t_1g_1g,...)$. The quotient of $EG$ by this right action is $BG$. The map from $EG$ to $BG = EG/G $ mapping an element of $EG$ to its right orbit is noted $p :EG \to BG$

Let $t_i$  be the map from $EG$ to $[0,1]$    sending $(t_0g_0,t_1g_1,...)$ to $t_i$;
 it induces also map $\underline {t}_i: BG \to (0,1]$ such that $p \circ t_i =\underline {t}_i$.
 
 Let $\omega \in H^1(BG, \tilde G)$ be the cohomology class 
 of the  $1$-cocycle defined as follows.
Let $(V_i)_{i= 0,1...}$ 
be the open cover 
$V_i ={t_i}^{-1}{(0,1]}$.
 Let $g_{ij}: V_i \cap V_j \to G$ be the map sending the class of $(t_0g_0,t_1g_1,\dots )$ to the element $g_ig_j^{-1}$; it is a $1$-cocycle whose cohomology class is $\omega$.

The infinite joint  construction of Milnor applied to $ G$ provides a canonical model for a universal principal $G$-bundle   $EG \to BG$.

{\bf 3.1.3. Numerable $G$-principal bundle.}
The classification theorem requires  the use of numerable covers. Such covers always exist on paracompact spaces. For this reason, we replace the category of topological spaces by the subcategory of paracompact spaces. The space $BG$ is also numerable.

An open cover $(V_i)_{i= 0,1...}$ of a topological space $X$ is numerable if there is a locally finite partition of the unity $f_i$ such that the closure of ${f_i}^{-1}(0,1]$ is contained in $V_i$. A $G$- principal bundle
$p: E \to X$ is numerable if it can defined by a $1$-cocycle over a numerable cover of $X$.

Two $G$-principal bundles  are numerably homotopic if they are connected by  a numerable homotopy.

For the technical points of  Milnor's joint construction of the classifying space, a good reference is the book of D. Husemoller: Fiber bundles, McGraw-Hill Series in higher Mathematics, 1966.

Recall that, 
for a toplogical space $X$,
then  $\tilde G$ is the sheaf of germs of continuous maps from $X$ to $G$.

\proclaim {Theorem 3.2} 

1) The map $p: EG \to BG$ is a $G$-principal bundle characterized by the element $\omega \in H^1(BG, \tilde G)$ defined at the end of 3.1.2.

2) This element is numerable.

3) Any numerable $G$-principal bundle 
with base space $X$   characterized by $\sigma \in H^1(X, \tilde G)$
is the pullback $ f^*\omega$ by a map $f: X \to BG$ whose homotopy class is well defined. Therefore:
$$H^1(X,\tilde G) = [ X, BG ],$$ 
where $[ X,  BG ]$ is the set of homotopy classes of continuous maps from a topological space  $X$ to the classifying space $BG$.

\endproclaim

\demo{Proof} For the proof, see [3], p. 142 - 143. We give some details of the proof of 3). 

A lemma asserts that, given a numerable cover $(U_i)_{i= 0,1,...}$, there is a locally finite countable partition $t_n$ of the unity, where $n\gneq 0$, such that each open $V_n= {t_n}^{-1}(01]$ is the disjoint union of open sets $V_{ni}$ contained in $U_i$.

Let $\sigma \in H^1(X,\tilde G)$ be a numerable $G$-principal bundle over $X$. Using the lemma, we can assume that $\sigma$ is the cohomology class of a $1$-cocycle $g_{mn}$ over the open cover 
$(U_n)_{n=0,1,...}$ defined by  $U_n = {t_n}^{-1}(01]$, where $t_n$ is a partition of the unity.

The required map $f: X \to BG$
  is defined by $f(x)$ equal to the $G$-orbit of  
$t_0g_{m0}(x), t_1g_{m1}(x),\dots$. Note that $g_{mn}(x)$ need to be defined only  for $x \in U_n$, because otherwise $t_n(x) =0$.$\square$
\enddemo

Before restating the theorem above to take care of base points,
we first consider the case of a  topological group $G$. The space $BG$ has a unique base point noted $*$ represented in the joint construction by the
infinite sequence consisting uniquely of the unit element of $G$.  
We have:
$$ \pi_i(BG,*) = \pi_{i-1}( G, 1), $$
where $1$ is the unit element of $G$. In particular, $\pi_1(BG,*)$ is the group consisting of the connected component of $G$. 

For instance, for a discrete group $\Gamma$, we have $\pi_1 (B\Gamma,*) = \pi_0(\Gamma,1) = \Gamma$; for $i >1$, $\pi_i (B\Gamma,*)$ is trivial. Therefore, $B\Gamma$ is an Eilenberg-MacLane complex $K(\Gamma,1)$. For more examples, see 3.5. below.

 We now come to the generalization of theorem 3.2, taking account of base points. In the notations above, the elements of $G$ are noted $g$ and the units (or the objects of $G$) are noted $y$.  
  A base point in $BG$ is determined as below by the choice of an element $*$ in the space of objects $Y$ of $G$. The choice  as base point of another element of $Y$ in the orbit of $*$ would lead to an isomorphic result.   Note that the choice of a base point which is not in the previous orbit could lead to a quite different result (see some of the examples in 3.5).
  
 Let $\tilde G$ be the sheaf of germs of continuous maps from $X$ to $G$, and consider a base point also noted $*$  in $X$. Then, $\tilde G_*$ will be the subsheaf of $\tilde G$ whose stalk above $*$ is the germ at $*$ of the constant map from $X$ to $* \in G$. The set 
 $H^1(X,\tilde G_*)$ is the subset of those elements of $H^1(X, \tilde G)$ represented by $1$-cocycles  $\phi_{ij}(u) \in \Cal Z^1(\Cal U,\tilde G)$ such that $\phi_{ij}(*)=1_*$. We denote by $[X,BG]_*$ the set of homotopy classes of  continuous maps from $X$ to $BG$ mapping $* \in X$ to $* \in BG$.

 The analogue of theorem 3.2, 3) above is:
 
 \proclaim{Theorem 3.3} Let $G_*$ be a topological groupoid with a base point $*$ and let $*$ be a base point in a topological space $X$. Then
$$H^1(X, \tilde G_*) =[X, BG]_*.  $$ 
\endproclaim

\demo {Proof} As in the proof of theorem 3.2, we only describe the map from $H^1(X,\tilde G_*)$ to $[X, BG]_*$. We choose a numerable cover
$(U_i)_{0,1,\dots}$ such that $* \in U_0$. Let $\sigma_* \in H^1(X , \tilde G_*)$ be a numerable $G$-principal bundle over the topological space $X$ and let $*$ be a base point in $X$. With the notations of the proof of theorem 3.2, $\sigma_*$ is the cohomology class of a $1$-cocycle $g_(mn)$ defined over the pointed open cover $(U_n)_{0,1,2,\dots}$, where $* \in U_0$.
The homotopy class of the  correspnding pointed map in $[X, BG]_*$ is represented by the map $f: X \to BG$ sending $x \in X$ to the $G$-orbit of $t_0g_{m0}(x), t_1g_{m1}(x), \dots$. 
$\square$
\enddemo

\proclaim{Corollary 3.4} The set of weak isomorphism classes of  homotopy classes of $(\Gamma \ltimes G)$-extensions of a space $X$ preserving a base point $* \in X$  is in bijection with the set $[X, BAut(\Gamma \ltimes G)]_*$.

In particular for $X= S^n$, this set is the group   $\pi_n(BAut(\Gamma \ltimes G),*)$.      
\endproclaim

\demo{Proof} According to part 1) of  theorem 2.11,  this set is equal to $H^1(X, \tilde {Aut}(\Gamma \ltimes G)_*)$. So we can apply theorem 3.3 to get the result. $\square$
\enddemo

 \noindent{\bf 3.5. Some examples.}

a) Assume that $G = G_1 \times G_2$ is the direct product  of two topological groups, then $BG = BG_1 \times BG_2$. The base point of $BG$ is the product of the base points of $BG_1$ and $BG_2$. We have then $\pi_i(BG,*) = \pi_i(BG_1,*) \times \pi_i(BG_2,*)$.

b) Consider the case of the  topological group $Aut(\Gamma \ltimes G)$ with base point $*$ the unit element. For a topological space $X$ with a base point $*$, then $H^1(X, \tilde{Aut}(\Gamma \ltimes G)_*))$ is the set of isomorphisms classes of pointed $(\Gamma \ltimes G)$-extensions of $X$. This set is in bijection with the set  $[X, BAut(\Gamma \ltimes G)]_*$. For $X = S^n$, we get:

$\pi_n(BAut(\Gamma \ltimes G),*)=\pi_{n-1}(Aut(\Gamma \ltimes G),*)$. Therefore  this is equal to $Aut(\rho)$ for $n=1$, to $\pi_{n-1}(G,*)$. As $G$ is simply connected, $\pi_i(G,*)=0$ for $i =1$ and $2$.

c) Let $G=\Gamma \ltimes T$ be the topological groupoid given by the action of a discrete goup $\Gamma$ acting by homeomorphisms of a topological space $T$. Then $BG = E\Gamma \times _\Gamma T$, the quotient of 
$E\Gamma \times T$ by the action of $\Gamma$ given by $\gamma.(e,t) =(e\gamma^{-1},\gamma t)$. This is a bundle with base space  $B\Gamma$ and fiber $T$. 

For instance, $B(\Gamma \ltimes G)$ is a bundle with base space $B\Gamma$ and fiber $G$. The exact homotopy sequence of this bundle gives the following result: $\pi_i(B(\Gamma \ltimes G),*)$ is equal to $\Gamma$ for $i=1$ and to $\pi_i(G,*)$ for $i>1$.

d) Consider an exact sequence  $1 \to K \to L \to M \to 1$  of topological groups. The fibration considered in 3.1.3  gives a map  $ BL \to BM$ which can be considered as the  projection of a fiber space with fiber   
$BK$. 

For instance consider the exact sequence  $1 \to G \to Aut(\rho) \ltimes G \to Aut(\rho) \to 1$; taking account of the isomorphism of  $Aut(\rho) \ltimes G$ with $Aut(\Gamma \ltimes G)$,   then 

\noindent$BAut(\Gamma  \ltimes G) \to BAut(\rho)$ is the projection of a fiber bundle with fiber $BG$.

e) Let us introduce base points in some of the examples c) above. Consider the topological groupoid $G = \Gamma \ltimes R^2$ associated  to the action of a discrete group $\Gamma$ acting effectively by rotations of $R^2$ keeping the origin fixed.

The map $\Gamma \ltimes R^2$ sending $(\gamma,x)$ to $\gamma$ induces a map $BG \to B\Gamma$ which is a homotopy equivalence because $R^2$ is contractible.

As base point $*$ of $BG$ let us first choose the origin $0 \in R^2$. 
Given a connected topological space $X$ with a base point $x_0\in X$, then the map $[(X,x_0),(BG,*)] \to [(X,x_0),(B\Gamma,1)]$ induced by the homotopy equivalence  $BG\to B\Gamma$ is a bijection.

Let us now choose  as base point $*$ a point of $R^2$ with positive norm, then the map described above with this new base point induces a map $BG \to B\Gamma$ which is a constant map. 

As an exercise, apply the same considerations to the case where $\Gamma \ltimes R^2$ is replaced by $SO(2) \ltimes R^2$, where the circle $SO(2)$ acts by rotations of $R^2$.

f) If  $G$ is a topological group, the base point is just the unit element $1$ of $G$. For $X= S^i$ with a base point $*$, we get:
$$ H^1(S^i, \tilde G_*)=\pi_i(BG,*).$$

\noindent {\bf 3.6. Crossed module.} We say that a group homomorphism $\mu :\Gamma \to E$ is a crossed module if an action of $E$ on $\Gamma$, noted $e.\gamma$ , is given such that :
$$ 1)\ \ \  e\mu(\gamma) e^{-1} = \mu (e.\gamma),$$ and
$$ 2) \ \ \ \mu(\gamma_1)\gamma_2 = \gamma_1\gamma_2 \gamma_1^{-1}.$$
This implies that the sequence  :
$$ 1 \to ker (\mu) \to \Gamma \to E \to coker (\mu) \to 1$$
is exact, where $coker ( \mu) = E/\mu(\Gamma)$.

In particular if $\mu$ is the inclusion of $\Gamma$ to a normal subgroup of $E$, then it is a crossed module for which the action  of $e$ on $\gamma$ is given by $e.\gamma =e\mu(\gamma)e^{-1}$.

If $E = Aut(\rho)$ and $\mu:\Gamma \to Aut(\rho)$ maps $\gamma$ to $ad(\gamma)$, then the exact sequence above is the one mentionned in 1.5.

\proclaim{3.7. Fundamental lemma} Let $p:\Cal G \to X$ be a
$(\Gamma \ltimes G)$-extension of a connected topological space $X$. Choose  a base point $*$ in $\Cal G$ projecting  to a point also noted $*$. It determines a base point $* \in B\Cal G$. We have an exact sequence:
$$ \pi_2(X,*) \to \Gamma \to \pi_1(B\Cal G,*) \to \pi_1(X,*) \to  1$$
which is mapped by a homomorphism $\phi$ to the exact sequence:
$$ 0 \to \Gamma_0 \to \Gamma \to Aut(\rho) \to Out(\rho) \to 1.$$
mapping $\pi_1(B\Cal G,*)$ to $Aut(\rho)$ and  $\gamma \in \Gamma$ to $ad(\gamma)\in Int(\rho)$.

There is a natural action of $\pi_1(B\Cal G,*)$ on the first sequence which can be considered as a crossed module.

The homomorphism $\phi$ is  a homomorphism of crossed modules.

 The group $\Gamma$ of the first exact sequence is mapped by the identity on the group $\Gamma$ of the second exact sequence. The map $\pi_2(X,*) \to \Gamma_0$ is  the composition of the map
$\pi_2(X,*) \to \Gamma$ with the identity of $\Gamma$, i.e. to the kernel of $\Gamma \to Aut(\rho)$.

\endproclaim

\demo{Proof} Recall that a pointed  $(\Gamma \times G)$-extension $p: \Cal G \to X $ is determined up to isomorphism by an element of the set $H^1(X, \tilde {Aut}(\Gamma \ltimes G)_*)$ (cf. theorem 2.11,1)). On the other hand, this set is equal to the set $[X, B\Cal G]_*$ of base point preserving homotopy classes of maps from $X$ to $B\Cal G$ (cf. theorem 3.3 and its corollary 3.4).

The map $p: \Cal G \to X$ induces a map still noted $ p :B\Cal G
 \to BX=X$ which can be considered as the projection of  a bundle over $X$ whose fiber over the base point $*$ is the classifying space  $B(\Gamma \ltimes G)$ of the $(\Gamma \ltimes G)$-extension of a point .  The first exact sequence is part of the homotopy exact sequence of this bundle, taking account that $\pi_1(B(\Gamma \ltimes G),*) =\Gamma$, because $G$ is  connected (see example c) in 3.5). The action of $\pi_1(B\Cal G)$ on this sequence is the usual action of the fundamental group of the total space of a fiber bundle  on the homotopy exact sequence of this bundle.
    
    The homomorphism $\phi$ to the second exact sequence is defined as follows. The homomorphisms of the homotopy groups of $B\Cal G$ are noted $\psi$ and those of the homotopy groups of $X$ are noted $\phi$.. The homomorphism $\psi :\pi_1(B \Cal G,*) \to Aut(\rho)$ is defined as follows: as $\pi_1(B\Cal G,*) = H^1(S^1, \tilde {Aut}(\Gamma \ltimes G)_*)= \pi_1(BAut(\Gamma \ltimes G),*) =\pi_0(Aut(\Gamma\ltimes G)) $: this is equal to $Aut(\rho)$ because $G$ is connected.

  The homomorphism $\phi: \pi_1(X,*) \to Out(\rho)$ is defined as follows. Given an element $\alpha \in \pi_1(X,*)$ represented by a loop $f : S^1 \to X$ based at $*$, 
  we can lift it to a loop $\tilde f$ in $B\Cal G$ based at $*$ representing an element $\tilde {\alpha} \in \pi_1(B\Cal G,*)$, 
 then $\phi(\alpha)$ 
  is the image of $\psi(\tilde{\alpha})$ by the projection $Aut(\rho) \to Out(\rho)$. Another lifting $\tilde f'$ would lead to an element $\alpha'$ differing from $\alpha$ by an element of the fundamental group of the fiber, i.e. by an element of $\Gamma$ which is in the kernel of $\phi$.
  $\square$

    \enddemo

\subheading {\bf  4. The homotopy type of the classifying spaces } 

In this section we prove the main results of this paper. Namely, we apply the preceeding considerations to the particular case where the topological groupoid is the groupoid  $\Gamma \ltimes G$ given by the action  of a dense subgroup $\Gamma$ of the Lie group $G$, acting on $G$ by left translations. 

In what follows, all the base points shall be noted by the same symbol $*$.

\proclaim{Theorem 4.1} 1) There exists a universal pointed $(\Gamma\ltimes G)$-extension of a pointed connected space

$$p^{\Gamma \ltimes G} : \ \  \Cal G_* ^{\Gamma \ltimes G} \to X_*^{\Gamma \ltimes G}$$ 
such that any pointed $(\Gamma \ltimes G)_*$-extension $\Cal G_* \to X_*$ is weakly isomorphic to the pullback  of this universal extension by a base point preserving continuous map $f: X_* \to X_*^{\Gamma \ltimes G}$.

2) The diagram of the fundamental lemma 3.7 where $p$ is replaced by $p^{\Gamma \ltimes G}$ is:
$$0\to \pi_2(X_*^{\Gamma \ltimes G},*)\to \Gamma \to \pi_1(B\Cal G_*^{\Gamma \ltimes G},*) \to \pi_1(X_*^{\Gamma \ltimes G},*) \to 1.$$ This exact sequence  is mapped isomorphically to the exact sequence:
$$0 \to \Gamma_0 \to \Gamma \to Aut(\rho) \to Out(\rho) \to 1.$$

3) The space $B\Cal G_*^{\Gamma \ltimes G}$ has the homotopy type of an Eilenberg-Maclane complex $K(Aut(\rho),1)$.

\endproclaim
\demo {Proof} 1) 
 The space  $BAut(\Gamma \ltimes G)_*$ is a classifying space for pointed $(\Gamma \ltimes G)$-extensions. Indeed by theorem 2.11.1), the set of pointed weak isomorphisms of $(\Gamma \ltimes G)$-extensions of a 
  topological space $X$ is equal to $H^1(X, \tilde {Aut}(\Gamma \ltimes G)_*)$, and this equal to $[X, BAut(\Gamma \ltimes G]]_*$ by theorem 3.3. Hence this space can be taken as $X_*^{\Gamma \ltimes G}$.
  
  Consider the following "exact sequence" of topological groupoids:
  $$ 1 \to (\Gamma \ltimes G) \to Aut(\Gamma \ltimes G) \ltimes (\Gamma \ltimes G) \to Aut(\Gamma \times G) \to 1,$$
  where the second arrow $p$ maps $(\alpha,(\gamma,g))$ to $\alpha$ 
  and the first arrow maps $(\gamma,g)$ to $(1,(\gamma,g))$.
   Passing to classifying spaces, we get a sequence:
  $$B(\Gamma \ltimes G) \to B( Aut(\Gamma \ltimes G) \ltimes (\Gamma\ltimes G)) 
 \to BAut(\Gamma \times G).$$
 The second arrow, still noted $p$, is a surjective map and the first one is injective. It is easy to check that $p$ is the projection of a fiber bundle with fiber $B(\Gamma \ltimes G)$. In the first sequence above, we can introduce base points by taking the unit element of  $Aut(\Gamma \ltimes G)$ and the element $(1 \ltimes 1) \in (\Gamma \ltimes G)$. This choice gives also base points in the second sequence. Therefore  $p=p^{\Gamma \ltimes G}: \Cal G ^{\Gamma \ltimes G}_* \to BAut(\Gamma \ltimes G)_*=X_*^{\Gamma \ltimes G}$.

2) Recall that the first sequence is part of the homotopy sequence
 of  the bundle $p^{\Gamma \ltimes G} : B\Cal G^{\Gamma \ltimes G} \to X^{\Gamma \ltimes G}$ with fiber $B(\Gamma \ltimes G)$,(see lemma 3.7).

    We first prove that $\psi: \pi_1 (B \Cal G^{\Gamma \ltimes G},*) \to Aut(\rho)$ is an isomorphism. Indeed $\pi_1 (B \Cal G^{\Gamma \ltimes G},* )=
H^1(S^1, \tilde {Aut}(\Gamma \ltimes G)_*)$ is equal to $Aut(\rho)$. This obviously implies that $\phi$ is an isomorphism.

On the other hand the homotopy exact sequence  of the fibration
$B\Cal G^{\Gamma \ltimes G} \to X^{\Gamma \ltimes G}_*$ gives:
$$ \pi_2(B(\Gamma \ltimes G),*) \to       \pi_2(B\Cal G^{\Gamma \ltimes G},*) \to \pi_2(X^{\Gamma \ltimes G},*) \to \Gamma. $$
The first group is $0$ (see example c) in 3.5), therefore the second  arrow is injective. This implies that $\pi_2(X^{\Gamma \ltimes G},*) = \Gamma_0$.

3) We have already proved in 2) that $\pi_1(B\Cal G^{\Gamma \ltimes G},*)=Aut(\rho)$; a similar argument tells us that  $\pi_2(B\Cal G^{\Gamma \ltimes G},*)=\pi_1(G,1)= 0$. We consider again the  homotopy exact sequence of the fiber bundle $B\Cal G^{\Gamma \ltimes G} \to  X^{\Gamma \ltimes G}$ with fiber $B(\Gamma \ltimes G)$.
To prove that all the other homotopy groups $\pi_i(B\Cal G^{\Gamma \ltimes G},*)$ for $i>2$  are equal to $0$, it will be sufficient to check that $\pi_{i+1}(X^{\Gamma \ltimes G},*) \to \pi_i(B(\Gamma \ltimes G),*)$ is an isomorphism for $i>1$.

We have $\pi_{i+1}(X^{\Gamma \ltimes G},*) = H^1(S^{i+1},\tilde{Aut}(\Gamma \ltimes G)_*) =[S^{i+1}, B(Aut(\Gamma \ltimes G)]_*   
=\pi_i(Aut(\Gamma \ltimes G),*) =\pi_i(G,*)$. On the other hand, we have proved that $\pi_i(B(\Gamma  \ltimes G),*) =\pi_i(G,*)$.
$\square$

\enddemo

\proclaim{Corollary 4.2} The homotopy groups $\pi_i(X^{\Gamma \ltimes G},*)$ are equal to $Out(\rho)$ for $i=1$, to $\Gamma_0$ for i=2, to $\pi_i(BG,*)$ for $i>2$. \endproclaim
\demo{Proof} This follows from the above proof.
\enddemo
Our next goal is the determination of the homotopy type of the classifying space $X^{\Gamma \ltimes G}$.

We have first to recall some basic notions.

\noindent {\bf 4.3. Postnikov decomposition}

An Eilenberg-MacLane complex $K(\pi,n)$ is a space such that
 $\pi_i(K(\pi,n))$ is equal to $\pi$ for $i=n$ and is trivial otherwise.
 It is connected for $n>0$. 
 For $n = 0$,  then  $K(\pi,0) $ is equal to the discrete space $\pi$ and  for $n>1$, then $\pi$ is abelian.
 
 Let $X$ be a connected space having the homotopy type of a CW-complex. 
 
 A Postnikov decomposition of $X$ is given by a sequence   of continuous maps of topological spaces:
 $$\dots  \to X_{(n)} \to X_{(n-1)} \to \dots \to X_{(2)} \to X_{(1)}$$
 
 and  continuous maps $f_n: X \to X_{(n)}$ such that:

 1) $X_{(1)} = K(\pi_1(X),1)$ and  for $ n>1$,
  $X_{(n)}$ is a topological space such that $\pi_p(X_{(n)})=0$ for $p>n$.
 
 2) The map $\alpha_{n} : X_{(n+1)} \to X_{(n)}$ is a fibration with fiber $K(\pi_{n+1}(X),n+1)$ for $n\geqq 1$.
 
 3) For $p�\leqq n$, the map $\pi_p(X) \to \pi_p(X_{(n)})$ induced by $f_n$ is 
 an isomorphism.
 
 4) The map $\alpha_n \circ f_{n+1}$ is homotopic to $f_n$.
 
 The choice of a base point $*$ in $X$ induces a base point in each $X_{(n)}$, and all the maps and homotopies are base point preserving.
 The group $\pi_1(X,*)$ acts on each $\pi_n(X,*)$. If this action is the identity, the Postnikov decomposition is called simple. The Postnikov decomposition of $B(Aut(\Gamma \times G))$ is not simple.

 Assume that the Postmikov decomposition of the connected topological space $X$  is simple. Then the fibration $\alpha_{n-1} : X_{(n)} \to X_{(n-1)}$ has a section if and only if an element  $p^{n} \in H^{n}(X,\pi_{n}(X))$ vanishes. Such an element is called a Postnikov invariant.

 In the next theorem, we describe  the Postnikov decomposition of $X^{\Gamma \ltimes G}$ and its homotopy type.

  \newpage

\proclaim{Theorem 4.4}

1) The second step $X_{(2)}^{\Gamma \ltimes G}$ of the Postnikov decomposition of $X^{\Gamma \ltimes G}$ is the total space of a fiber bundle with base space $X_{(1)}^{\Gamma \ltimes G}=K(Out(\rho),1)$ and fiber $K(\Gamma_0,2)$. It is the classifying space of the crossed module $0 \to \Gamma_0 \to \Gamma \to Aut(\rho) \to Out(\rho)\to 0$,

2) The homotopic fiber  of the projection of $X^{\Gamma \ltimes G}$ to the second step of its Postnikov decomposition is homotopically equivalent to $BG$. In fact, $X_{(3)}^{\Gamma \ltimes G}=X_{(2)}^{\Gamma \ltimes G}$.
\endproclaim
\demo{Proof} 1) The first sentence follows from the definition. For the second one, see J.-L. Loday: Spaces with finitely many non-trivial homotopy groups, J. Pure Applied Algebra, 24, (1982), 179-202.

2) By definition, $X_{(n)}^{\Gamma \ltimes G}$ is $2$-connected for $n>2$. The  projection  $f_2 : X^{\Gamma \ltimes G} \to X_{(2)}^{\Gamma \ltimes G}$ is homotopiquement equivalent to a fibration still noted $f_2$ (see J.-P. Serre: Homologie singuli\`ere des espaces fibr\'es, Ann. of Math., vol.54, 1951, pp. 425-503, Chapitre IV).

We have to prove that the fiber $F$ has the homotopy type of $BG$. First, mention that $\pi_2(G,1) = 0$  for any simply connected Lie group $G$. Hence, $\pi_i(BG,*)$ is trivial for $i<
3$. This implies that $X_{(3)}^{\Gamma \ltimes G}=X_{(2)}^{\Gamma \ltimes G}$.

The exact sequence of the fibration $f_2$ implies
that the inclusion of $F$ in $X^{\Gamma \ltimes G}$ induces an isomorphism  from the homotopy groups of $F$ and those of $BG$.$\square$
\enddemo

\proclaim{Corollary 4.5} If $G$ is contractible, then $X^{\Gamma \ltimes G} =X_{(2)}^{\Gamma \ltimes G}$. \endproclaim

\newpage

\heading{III.  Examples.}\endheading

\noindent{\bf 1.The case where the center $\Gamma_0$ of $\Gamma$ is discrete in $G$.}

We first prove a lemma which shall be useful in this section.

\proclaim{1.1. Lemma} Consider a sequence of groups $K_0 \to K \to A$, where $\iota : K_0 \to K$ and $\mu: K\to A$ are injective homomorphisms of groups. We assume that $\iota(K_0)$ is a normal subgroup of $K$, that $\mu(K)$ is a normal subgroup of $A$ and that $(\mu \circ  \iota) (K_0)$ is also a normal subgroup of $A$. The above sequence of groups induces the exact sequence of groups :
$$ 1\to K/K_0 \to A/K_0 \to A/K \to 1$$ so that  $A/K = A/K_0\ \  /\ \ K/K_0$. \endproclaim

\demo{Proof} We have the  sequence of groups:
$$1\to K/K_0 \to A/K_0 \to A/K \to 1,$$
where the first homomorphism $\overline \mu$ is defined by $\overline \mu(kK_0)= \mu(k)K_0$ and the second one $\overline \iota$ by $\overline \iota (aK_0) =aK$. We check the exactness of the sequence. The map $\overline \mu$ is injective because $\mu$ is injective; the element $aK_0 \in A/K_0$ is in the kernel of $\overline i$ iff $a \in K$; in that case it is the image by $\overline \mu$ of $a K_0$. It is clear the $\overline i$ is surjective because $aK \in A/K$ is the image of 
$aK_0$. $\square$

Let $\Gamma_0$ be the center of $\Gamma$. We consider the case where $\rho(\Gamma_0)$ is a dicrete subgroup of $G$. We have proved in II,1.7  that the homomorphism $\mu: \Gamma \to Aut(\Gamma\ltimes G) =Aut(\rho) \ltimes G$ defined by $\mu(\gamma)=(ad(\gamma), \rho(\gamma^{-1}))$ maps $\Gamma$ to a normal subgroup $\mu(\Gamma)$ of $Aut(\Gamma \ltimes G)$. Note that if the center of $G$ is discrete, then $\rho(\Gamma_0)$ is discrete in $G$. The converse is not true in general (see the examples below).

\proclaim{1.2. Theorem} If $\rho(\Gamma_0)$ is a discrete subgroup of $G$, then $\mu(\Gamma)$ is a discrete subgroup of $Aut(\Gamma \ltimes G)$. 

Therefore $Aut(\Gamma \ltimes G)/\mu(\Gamma) = Out(\Gamma \ltimes G)$ is a Lie group $H$. 
\endproclaim

\demo{Proof} Consider the composition of $\mu$ with the projection to the first factor $Aut(\rho)$. Its image is 
$Ad(\Gamma)$ and its kernel is $\Gamma_0$. Therefore, the intersection of  $\mu(\Gamma)$ with the connected component $G$ of the identity of the Lie group $Aut(\Gamma \ltimes G)$ is equal to $\rho(\Gamma_0)$. So, if this group is discrete, then $\mu(\Gamma)$ is discrete. $\square$\enddemo

The main purpose of this section is to describe the structure of the Lie group $H$ (which is not connected in general).
The group $Out(\rho)$ acts on $G/\rho(\Gamma_0)$ as follows. Given $\tilde \alpha \in Out(\rho)$ which is the image of $\alpha \in Aut(\rho)$, its action on $G/\rho(\Gamma_0)$ is equal to $\tilde \alpha(g\rho(\Gamma_0)) = \overline \alpha(g) \rho(\Gamma_0)$. The action of $Aut(\rho)$ on $G$ leaves $\rho(\Gamma_0)$ invariant  because $Aut(\rho)$ acts by automorphisms on $\rho(\Gamma)$

\proclaim{1.3. Theorem } 
We have the exact sequence of groups:
$$1 \to \Gamma /\Gamma_0 \to Aut(\Gamma \ltimes G)/\mu(\Gamma_0) \to Aut(\Gamma \ltimes G)/\mu(\Gamma) \to 1,.$$ where the first homomorphism is noted  $i$ and the second one $p$.

This sequence is isomorphic to the exact sequence of groups: $$1 \to Int(\rho) \to Aut(\rho)\ltimes G/\rho(\Gamma_0)  \to Out/\rho)\ltimes G/\rho(\Gamma_o) \to 1, $$
where the first homomorphism maps $Int(\rho)$ to the subgroup $Int(\rho)  \ltimes 1$

The group $H$ is an extension of the connected Lie group $H_0=G/\rho(\Gamma_0)$ by the discrete group $Out(\rho)$ and it is isomorphic to $Out(\rho) \ltimes G/\rho(\Gamma_0)$.

 \endproclaim 

\demo{Proof} The group $\mu(\Gamma)$ is an invariant subgroup of $Aut(\Gamma \ltimes G)$. So we can apply  lemma 1.1
to the sequence of groups:
$$\Gamma_0 \to \Gamma \to Aut(\Gamma \ltimes G),$$
and we get the exact sequence :
$$1 \to\Gamma/\Gamma_0 \to Aut(\Gamma\ltimes G)/\mu(\Gamma_0) \to Aut(\Gamma \ltimes G)/\mu(\Gamma)\to 1.$$

The group $\Gamma/\Gamma_0$ is isomorphic to $Int(\rho)$. Indeed we have the exact sequence: $1 \to \Gamma_0 \to \Gamma\to Aut(\rho) \to Out(\rho) \to 1$, where the homomorphism $\Gamma \to Aut(\rho)$ maps $\gamma$ to $ad(\gamma)$ (see I, 1.5).

Recall that $Aut(\Gamma \ltimes G)$ is isomorphic to $Aut(\rho) \ltimes G$. We also have  $\mu(\gamma) =(ad(\gamma),\rho(\gamma)^{-1})$. Hence   $\mu(\Gamma_0) = (1,\rho(\Gamma_0))$. 

 Therefore $Aut(\Gamma \ltimes G)/\mu(\Gamma_0) = Aut(\rho)\ltimes G/\rho(\Gamma_0)$. The image of the homomorphism $\Gamma /\Gamma_0  =Int(\rho)\to Aut(\Gamma \ltimes G)/\mu(\Gamma_0) $ maps $\alpha \in Int(\rho)$  to  $ \alpha \ltimes 1$, where $1$ is the unit element of the group $G/\rho(\Gamma_0)$ . 
    As $Aut(\rho) / Int(\rho) = Out(\rho)$ we have $H = Aut(\Gamma \ltimes G)/\mu(\Gamma) =Out(\rho) \ltimes G/\rho(\Gamma_0)$. $\square$ \enddemo

\proclaim{1.4. Corollary} If $\rho(\Gamma_0)$ is discrete in $G$, then $BOut(\Gamma \ltimes G)$ is isomorphic to the classifying space $BH$ of the group $H = Aut(\Gamma \ltimes G)/\mu(\Gamma)$.
Therefore, there is a bijection between the set  of homotopy classes  of equivalence classes of $(\Gamma \ltimes G)$-extensions of a topological space $X$ and the set of homotopy classes of continuous maps from $X$ to $BH$.

\endproclaim

\proclaim{1.5. Theorem} 1) The homotopy groups $\pi_i(BOut(\Gamma \ltimes G),1)$ are equal to $Out(\rho)$ for $i=1$, to $\Gamma_0$ for $i=2$ and to $\pi_{i-1}(G,1)$ for $i>2$.

2) The second step $BOut(\Gamma \times G)_{(2)})$ of the Postnikov decomposition of $BOut(\Gamma \ltimes G)$ is equal to 
$K(Out(\rho),1) \times K(\Gamma_0,2)$. The homotopic fiber of the projection $BOut(\Gamma \ltimes G) \to Bout(\Gamma \ltimes G)_{(2)}$ is homotopically equivalent to $BG$.

The natural projection $BAut(\Gamma \ltimes G) \to BOut(\Gamma \ltimes G)$ with kernel $B\Gamma$ induces a map of their Postnikov decomposition.

\endproclaim

\demo{Proof} 1) The exact sequence of theorem 1.3 induces a sequence of classifying spaces:
$$ B(\Gamma/\Gamma_0) \to B(Aut(\rho)\ltimes G/\rho(\Gamma_0)) \to BOut(\Gamma \ltimes G)$$
where the second arrow is the projection of a fiber bundle with fiber 
$B(\Gamma / \Gamma_0)$.
To get the result, write the homotopy exact sequence of this bundle, taking account of the equality $\Gamma /\Gamma_0 = Int(\rho)$.

\noindent 2) The exact sequence of groups  $ \Gamma \to Aut(\Gamma \ltimes G) \to Out(\Gamma \ltimes G) \to 1$, where the first arrow is $\mu$, induces a sequence $B\Gamma \to BAut(\Gamma \ltimes G) \to BOut(\Gamma \ltimes G)$, where the second arrow is the projection of a fiber bundle with fiber $B\Gamma$. In turn, this induces maps from the Postnikov decompositions of those spaces.

The first step is equal to the sequence: $K(\Gamma,1) \to K(Aut(\rho),1) \to K(Out(\rho),1)$.
The second step is the following: $BAut(\Gamma \ltimes G)_{(2)}$ is a fiber bundle with fiber $K(\Gamma_0,2)$ and base space $K(Aut(\rho),1)$; and $BOut(\Gamma \ltimes G)_{(2)}$ is a fiber bundle with fiber $K(\Gamma_0,2)$ and base space $K(Out(\rho),1)$. The projection $BAut\Gamma \times G) \to BOut(\Gamma \ltimes G)$ induces a bundle  map which is an isomorphism
on the fibers $K(\Gamma_0,2)$. As $Out(\rho)$ acts by the identity on $\Gamma_0$, this implies that $B(Out(\Gamma \ltimes G)_{(2)} = K(Out(\rho),1) \times K(\Gamma_0,2)$.

In both cases, the fiber of the projection to the second step is isomorphic to $BG$, but, in the case of $BOut(\Gamma \ltimes G)$, as $Out(\rho)$ acts trivially, its Postnikov decomposition is simple.$\square$

\enddemo
\proclaim{1.6. Proposition} If the center of $\Gamma$ is discrete, if $\Gamma$ is countable and if $X$ is a connected complete paracompact smooth manifold, then any $(\Gamma \ltimes G)$-extension of $X$ is equivalent  to the holonomy pseudogroup of a Riemannian foliation on a complete Riemannian manifold. \endproclaim 

\demo{Proof} Let $p: E \to X$ be a principal $H$-bundle with basis $X$. Our first aim is to construct on $E$ a complete 
Riemannian metric invariant by the right action of $H$. 
 Let $g$ be a complete Riemannian metric on $X$ and let $g_H$ be a Riemannian metric on $H$ invariant by the right translations.
Let $\Cal U$ be a cover of $X$ by open sets $U_i$ such that $p^{-1}(U_i)$ is isomorphic to a product  $U_i \times H$. 
Let $\lambda_i$ be partition of the unity subordinated to the open cover $\Cal U$. Then, the desired metric is equal to $\sum \lambda_i g_i \times g_H$, where $g_i$ is the restriction of $g$ to $U_i$.

Let $E_0$ be a connected component of $E$ and let $Out(\rho)_0$ be the subgroup of $Out(\rho)$ leaving $E_0$ invariant; this subgroup is countable, because it is a quotient of $\Gamma$. Let $S$ be a $2$-dimensional  surface with a complete Riemannian metric $g_S$ such that there is a surjective homomorphism $\phi: \pi_1(X) \to Out(\rho_0)$.
Let $\hat S$ be its universal cover with its induced metric. On the manifold $\hat S \times E_0$ with the product Riemannian metric, the horizontal foliation is invariant by the action of $\pi_1(S)$ acting on $\hat S$ by covering translations and on $E_0$ by right translations via $\phi$. The quotient foliation is the Riemannian foliation we are looking for.

\enddemo

\noindent {\bf 1.7. Examples where $\rho(\Gamma_0)$ is discrete in $G$.}

a) The center of any simply connected compact Lie group $G$ is finite, hence discrete. So for  any dense subgroup $\Gamma$ of $G$, then $\mu(\Gamma)$ is discrete in $Aut(\Gamma \ltimes G)$.

b) It may happen that $\rho(\Gamma_0)$ is discrete in a simply connected Lie group $G$ whose center is not discrete. The following is an example where $\Gamma_0$ is trivial. Let $G = K \times R$, where $K$ is a simply connected compact Lie group, for instance  $SU(2)$.  It is known (see P. de la Harpe, Free groups in linear groups, Ens. Mat. 29 (1963), pp 129 - 144) that there exist  a finitely generated subgroup (even a free subgroup) $\Gamma$ dense in $K$. Let $r: \Gamma \to K$ be the inclusion of $ \Gamma$ in $K$, and let 
 $h : \Gamma \to R$ be a homomorphism whose image is dense in $R$.

Let $\rho(\Gamma)$ be the subgroup of $K \times R$ which is the image of $\Gamma$ by the homomorphism $\rho:\gamma \mapsto (r(\gamma),h(\gamma))$. 
Then $\rho(\Gamma)$ is dense in $K\times R$. Indeed,  the closure of $\rho(\Gamma)$ is a closed Lie subgroup $H$ of $K \times R$. If $ dim H= dim K$, then the kernel of the   projection on $K$ would be equal to a a non trivial finite subgroup of $
1 \times R$ which does not exist. Therefore, $dim H > dimK$, i.e $H =K\times R$.

Hence $Aut(\rho)=Aut(r) \times Aut(h(r))$.

c) Consider the case where the simply connected Lie group $G$ is the group $Spin(n)$, the 2-fold covering of $SO(n)$. The center of $Spin(n)$ is equal to $Z_2$ for $n=2k+1$, to $Z_4$ for $n=4k+2$ and to $Z_2 +Z_2$  for $n =4k$.

Let $\rho$ be the inclusion of a group $ \Gamma $ in the quotient of $Spin(n)$ by its center. Then, the inverse image of $\rho(\Gamma)$ by the projection of $Spin(n)$ in this quotient. Then, $p^{-1}(\rho(\Gamma))$ is isomorphic to the product of $\Gamma$ by the center of $Spin(n)$.

\noindent {\bf 2. The case where $G$ is abelian.}

Any abelian simply connected Lie group
 is isomorphic to $R^n$. The group of automorphisms of $R^n$ is equal to $Gl(n,R)$.

\noindent{\bf 2.1.}  We consider first the case $n=1$. The group $Gl(1,R)$ is the group $R^*$ of non-zero real numbers acting on $R$ by homotheties. The group $Aut(\rho)$ contains always the subgroup of order $2$ consisting of the identity and the central symmetry. We shall assume that $\Gamma$  contains, at least, an irrational number $\alpha$, so it contains the subgroup $Z \oplus Z\alpha$. We shall also assume that $\Gamma$ is of finite rank.

\proclaim{Proposition 2.2}
a) Assume that  the subgroup $\Gamma$ of $R$ contains  a transcendental number 
 $\alpha$, then $Aut(\rho) =\{1,-1\}$.

b) Assume that $\Gamma$ is finitely generated and let $N$ be the rank of $\Gamma$. Let  us call $\Lambda$  the group $Aut(\rho)$. Let $\kappa$ be the subfield of the field $R$ generated by $\Lambda$. Then, $\kappa$ is a number field whose degree is a divisor of $N$ and $\Lambda$ is a subgroup of the group of the units of the ring of integers of $\kappa$. 
 \endproclaim
 
 \demo{Proof}
a) Assume that $Aut(\rho)$ contains an element $\lambda$ which is not the identity and the central symmetry. Then we would have $\lambda^n(\alpha) \neq \pm 1$ for all $n>0$, contradicting the assumption
that the rank of $\Gamma$ is finite.

b) Let $\lambda \in \Lambda$; it can be considered either as an automorphsm of $R$ or as an automorphism of the free $Z$-module $\Gamma$. It is,then, a root of its characteristic polynomial, which is of degree $N$ with integral coefficients, the constant term being $1$ or $-1$. So the real number $\lambda$ is a unit of  the field $Q(\lambda)$.

Let $\kappa$ be the subfield of the field $Q(\Lambda)$ of $R$ generated by $\Lambda$. Its elements act by multiplication on the vector subspace $\Gamma \otimes Q$  of dimension $N$ over $Q$. Hence,$\Gamma \otimes Q$ can also be considered as a vector over $\kappa$. The product of its dimension over $\kappa$ by the degree $r$ of $\kappa$ over $Q$ is its dimension $N$ over $Q$; therefore, $r$ is a divisor of $N$.
$\square$

\enddemo

As a reference for basic notions of number theory, see for instance:  P. Samuel, Th\'eorie alg\'ebrique des nombres, Hermann Paris, (1971).

\proclaim{2.3. Corollary} If the rank $N$ of $ \Gamma$ is a prime number, then $Aut(\rho)$ is either equal to $\{+1,-1\}$ or $\Gamma$ is a multiple of the submodule  of the module of integers of a number field of degree $N$. So it can be considered as a vector space over $\kappa$. The product of its dimension over $\kappa$ by the degree $r$ of $\kappa$ is of dimension $N$ over $Q$, so that $r$ is a divisor of $N$.

\endproclaim

\demo{Proof} If $ \Lambda = Aut(\rho)$  contains an element $\lambda \neq \pm 1$, the $\lambda$ is a unit of a number field of degree prime to $N$. Let $\gamma \neq 0$ be an element of $\Gamma$; then the elements $\gamma,  \lambda \gamma,\dots, \lambda ^{N-1} \gamma$ are linearly independent in $Q$, hence they form a basis of $\Gamma \otimes Q$. Let $\Cal O$ bthe ring of integers of the field $\kappa$ above; it is a module of rank $N$ over $Z$, and $\Cal O \otimes Q = \kappa$. As $\gamma^{-1}\Gamma \subset \kappa$, there is a suitable integer $r$ such that $r \gamma^{-1}\Gamma \subset \Cal O$. $\square$

\enddemo
 
Note that in general for an integer $N$, then $Aut(\rho) \neq \pm 1$ only if $\Gamma$ is a sum of multiples of $Z$ modules  (not all monogenes) generated  by multiplicative subgroups of the group of units of the number field $\kappa$.

\noindent{\bf 2.4. The case where $G=R^2$ or  C.}

We can identify $R^2$ to $C$ by the map $(x,y) \to x+iy$.

Assume that $\Gamma \subset C$ is generated by  $Z \times Zi$ 
and by $(1,\alpha i)$, where $\alpha$ is an irrational number $<1$.
Then $\Gamma$ is dense in $C$, because in the torus quotient  of $C$ the line through the origine of slope $\alpha$ projects to a dense geodesic. Let $\rho: \Gamma \to C$ be the inclusion of $\Gamma$ in $C$.

Let $p: x+yi \mapsto y$ be the projection $C $ to its imaginary part. The image of $\Gamma$ by $p$ is generated by $1$ and $\alpha$. Let $\rho_0$ be the inclusion of $p(\Gamma)$ in $R$. The kernel of $p$ restricted to $\Gamma$ is isomorphic to $Z$, because it is generated by the translation $(x,y) \to (x+1,y)$. The projection $p$ induces a surjective homomorphism $Aut(\rho) \to Aut(\rho_0)$, Therefore, we have the exact sequence $1 \to Z \to Aut(\rho) \to Aut(\rho_0) \to 1$. In fact, $Aut(\rho)$ is isomorphic to $Z \times Aut(\rho_0)$: it is clear that $Z$ commutes with $Aut(\rho_0)$; on the other hand, the map $\lambda: Aut(\rho_0) \to Aut(\rho_0)$ generated by $\lambda(\alpha)= (0,\alpha)$ and $\lambda(1) = (0,1)$ induces a lifting of $p$.

For instance, if $\alpha$ is transcendental, then $Aut(\rho_0) = Z_2$; hence $Aut(\rho) =Z \times Z_2$. If $\alpha = \sqrt {2}$, then $Aut(\rho_0)= Z \oplus Z_2$ (see II, 3.11), hence $Aut(\rho)= Z \times(Z\oplus Z_2)$.

\noindent
{\bf 3. The case where $G$ is nilpotent.} 

In this section we assume that G is a simply connected nilpotent Lie group and that the dense subgroup $\Gamma$ is finitely generated (see [4], 4.2). Note that the center $\Gamma_0$ of $\Gamma$ is never discrete in $G$ and that $\Gamma$ is a nilpotent group without torsion.

\proclaim{3.1. Mal'cev theorem} Let $\Gamma$ be a finitely generated dense subgroup of a simply connected nilpotent Lie group $G$. Then there is a simply connected nilpotent Lie group $\tilde G$ containing a cocompact discrete subgroup $\tilde \Gamma$ and a surjective homomorphism $q_0:  \tilde G \to G$ whose restriction to $\tilde \Gamma$ is an isomorphism  $\tilde \Gamma \to \Gamma$.
\endproclaim 

\demo {Proof}See A.I. Mal'cev, Nilpotent torsion free groups, Izvestiya Akad. Nauk SSSR, Ser. Mat. 13. (1949), 201 - 212). $\square$
\enddemo
 
\proclaim{3.2. Corollary}For every dense finitely generated subgroup $\Gamma$ of a simply connected nilpotent Lie group $G$, then $\Gamma \ltimes G$ is the holonomy groupoid of a foliation on a compact manifold. \endproclaim 

\demo{Proof} Consider the foliation on $\tilde G$ whose leaves are the fibers of the submersion $q_0$. The action  of $\tilde {\Gamma}$ maps leaves to leaves. The quotient of this foliation
 by the cocompact subgroup $\tilde {\Gamma}$ is a compact manifold carrying the quotient foliation. Its holonomy groupoid is $\Gamma \ltimes G$. 
 $\square$ \enddemo

  The inclusion of $\tilde \Gamma $ in $\tilde G$ is noted $\tilde\rho$. Note that the homomorphism $q_0$ is a submersion whose kernel is a simply connected nilpotent Lie group.

 Moreover, there is an isomorphism from $Aut(\rho)$ to a group of automorphisms $Aut(\tilde{\rho})$ of $\tilde \Gamma$ mapping $\alpha \in Aut(\rho)$ to an automorphism $\tilde \alpha$, such that $q_0 \circ \tilde \alpha = \alpha $. 
  We also have the exact sequence:
$$1 \to Int(\tilde{\rho}) \to Aut(\tilde{\rho}) \to Out(\tilde{\rho}) \to 1,$$ projecting isomorphically to the exact sequence in II, 1.5.

Note that an abelian group is nilpotent. If $G=R$ and the dense subgroup $\Gamma$ has a basis $x_1,\dots, x_n$, then $\tilde G$ can be chosen to be $R^n$, and $q_0$ can be defined by $q_0(e_i) =x_i$, where $e_i$ is a basis of $R^n$. Then, $\tilde \Gamma$ can be identified to $Z^n$.

The group $Out(\tilde{\rho})$ acts on $\tilde G/\tilde{\rho}(\tilde \Gamma)$ as follows: an element of $Out(\tilde{\rho})$ is the image of an element $\tilde {\alpha} \in Aut(\tilde{\rho})$; its action on the coset $\tilde g \rho(\tilde \Gamma)$ is equal to the coset $\tilde \alpha (\tilde g)\rho(\tilde \Gamma)$.

We can apply  theorem 1.3 to the group $\tilde G$, because the subgroup $\tilde \Gamma$  is discrete in $\tilde G$ and a fortiori its center $\tilde \Gamma_0$. Therefore, we get the following 
 result:

\proclaim{3.3. Theorem}
Let $H$ be the Lie group quotient of $ Aut(\tilde{\rho}) \ltimes \tilde G$ by the discrete normal subgroup $\tilde\mu(\tilde\Gamma)$, where $\tilde \mu( \tilde \gamma) =(ad(\tilde\gamma), \tilde \gamma^{-1})$. The group $H$ is isomorphic to the semi-direct product  $Out(\tilde \rho)\ltimes \tilde G/\rho(\tilde \Gamma_0)$, where the discrete group $Out(\tilde{\rho})$ acts on $H_0 =\tilde G /\tilde\Gamma_0$ as described above.

\endproclaim
\proclaim {3.4. Theorem }
The Lie group $H$ carries a foliation $\Cal F$ invariant  by the  its left translations. Its holonomy groupoid is isomorphic to the groupoid given by the action of $\mu( \Gamma)$  on  $Aut(\rho) \ltimes  G$.
\endproclaim

\demo{Proof} Let $\Cal F_0$ be the foliation on $\tilde G$ whose leaves are the fibers of the submersion $q_0: \tilde G \to G$. This foliation is invariant by $\tilde G$ acting ont itself by left translations. It can be extended to a foliation $\Cal F$ on $Aut(\tilde \rho) \ltimes \tilde G)$ invariant by left translations. The space of leaves of $\Cal F$ is isomorphic to $Aut(\rho) \ltimes G$. The foliation $\Cal F$ is also invariant by the sugroup $\tilde \mu \tilde \Gamma)$. Therefore, we get on $H$ a foliation $\tilde \Cal F$ left invariant by $H$. Its holonomy groupoid is isomorphic to the action of $\mu(\Gamma)$  on $Aut(\rho) \ltimes G$. $\square$

\enddemo

\proclaim{3.5. Corollary} The classifying space $BH$ has the homotopy type of $BOut(\Gamma \ltimes G)$.

\endproclaim

\demo{Proof} The proof is similar to the proof of corollary  1.4. We first note that the kernel of the homomorphism $q_0 : \tilde G \to G$ is contractible, because it is a connected nilpotent Lie subgroup. This implies that the same is true for $q : Aut(\tilde \rho) \ltimes \tilde G \to Aut(\rho) \ltimes G$. Therefore, the leaves of the foliation $ \Cal F$ on $H$ have contractible holonomy cover.

 According to [4], 3.2.4,    $BH$ has the homotopy type of
 $BOut(\tilde\Gamma \ltimes \tilde G)$. This space is homotopically equivalent to $BOut(\Gamma\ltimes G)$.
$\square$
\enddemo

\noindent {\bf 4. The case where $G$ is solvable.}

Consider  the simply connected solvable Lie group $GA$ of affine transformations of $R$ preserving the orientation. It can be described as the group of matrices of the form:

$$
\left (\matrix
 a & b\\
 0 & 1
 \endmatrix \right)
 $$
where $a>0$. The map $GA \to {R_+}^* \ltimes R$ associating to such a matrix    the element $(a,b) $ of  the semi-direct group ${R_+}^* \ltimes R$, where the group of ${R_+}^*$ acts on $R$
 by positive homotheties, is an isomorphism. 
 
 In the appendix E of Etienne Ghys in [6], pp. 297-314, it is proved using number theory considerations that $GA$ contains a dense subgroup $\Gamma_n$, for each $n>1$, which is isomorphic to 
 $Z^{n-1} \ltimes Z^n$.

Following Yves Carri\`ere, we consider  the cases where
$n$ is equal to $2$ and $3$ (see the appendix A in [6], Example 2.1, p. 221).

Let $A$ be a matrix in $Sl_2(Z)$ with  positive trace. It has two positive eigenvalues $0 <\lambda_1 <1<\lambda_2$. Moreover we have $\lambda_1\lambda_2 =1$. Let $V_1=(1,\lambda_1)$ and $V_2=(1,\lambda_2)$ be the corresponding eigenvectors in $R^2$; therefore we have $A(V_1) =\lambda_1V_1$ and $A(V_2) = \lambda_2 V_2$.
 
 The straight lines parallel to $V_1$ and $V_2$ give on the quotient $T^2 = R^2/Z^2$ two foliations $F_1$ and $F_2$ with dense leaves   because $\lambda_1$ and $\lambda_2$ are irrational. The automorphism $A: R^2 \to R^2$ induces  automorphisms of $F_1$ and $ F_2$.

As a concrete example, consider the case where $A$ the matrix where $a=2$ and $b=1$. Then the eigenvalues of $A$ are equal to $3 \pm
\sqrt 5 /2$.

The subgroups of $GA$ generated by the elements 
$(1,\lambda_1)$ and $(1, \lambda_2)$ in $R_+^* �\times R$ are isomorphic to $Z$ and 
are invariant by the action of $A$ by homothetie.

The dense subgroup $\Gamma_2$ in $GA$ is isomorphic to $Z \ltimes  Z^2$, where the first factor is generated by $A$ and the second one by $V_1$ and $V_2$.

The dense subgroup $\Gamma_3$ of $GA$ is isomorphic to $Z^2 \ltimes Z^3$. The first generator of $Z^2$ is noted $A$ and the second one $A'$. The generators of $ Z^3$ are noted $V_1$, $V_2$ and $V_3$. The action of $A$ on $V_1$ and $V_2$ is given by $A(V_1) = \lambda_1V_1$ and $A(V_2)=\lambda_2V_2$ and is the identity on $V_3$. The action of $A'$ on $V_1$ and $V_2$ is the identity and $A'(V_3) = A^{-1}(V_3)$.

For more details, see the appendix A of Yves Carri\`ere (cf. loc. cit.).
The general case is explained in the appendix E of Etienne Ghys.
One can find many interesting examples of subgroups of $GA$ in the paper of Gael Meigniez, pp.109-146, in the book "Integrable Systems and Foliations" in Progress of Mathematics, Vol. 145, edited by C.Albert, R. Broouzet, J.-P. Dufour.


\newpage

{\bf Postface.}

In this postface we give a brief historical account on the theory of Riemannian foliations and of the theory of pseudogroups of local isometries.

At the end we evoke vaguely the wealth of problems that could be studied and which involve  Lie groups, number theory, algebraic topology, complex geometry, analysis, etc..

 The notion of Riemannian foliation was introduced by Bruce Reinhart in [7], under the name of foliation with bundle like metric. He gave several characterisations of such  foliations $\Cal F$ on a compact Riemmanian manifold $M$. One of them is the following: a geodesic orthogonal to a leaf of $\Cal F$ at one point remains orthogonal to any other leaf it meets. He also proved,  among other properties, that all the leaves have isomorphic simply connected covers if $M$ is connected; and, as well, that the space $X$ of the closures of the leaves with the quotient topology is a Hausdorff space. 
  
 In the late seventies, Molino proved
 a  deep structure theorem for Riemannian foliations on complete Riemannian manifolds $M$ (see [5],[6]). In particular,  he associated to each such foliations on a connected manifold a Lie algebra. The translation of Molino theory to our  point of view,  in terms of pseudogroups of local isometries, was well explained in [8].
 
 With the notations above let $p: M \to X$ be the map mapping a point of $M$ to the closure of the leaf containing this point. The space $X$ has a natural stratification; roughly speaking, a stratum consisting of all closures of leaves having the same type. Given a stratum $X_i$ in $X$, a fundamental problem is to describe the holonomy  pseudogroup of the restriction $\Cal F_i$ of $\Cal F$ to $M_i =p^{-1}(X_i)$. This pseudogroup can be considered as a locally trivial extension of the topological space $X_i$ by the holonomy pseudogroup of the restriction of  $\Cal F_i$ to the closure of one leaf. One should also understand how the strata fit together. 
 
 We now come back to our program of studying directly the holonomy pseudogroups of Riemannian foliations (see [4]). The holonomy pseudogroup $P$ of a Riemannian foliation $\Cal F$ on a Riemannian manifold $M$
can be considered as a pseudogroup of local isometries of a Riemannian
 manifold, namely a transversal submanifold $X$ of dimension equal to the codimension $n$ of $\Cal F$ and meeting all the leaves; the submanifold $X$ is in general not connected and is endowed with the transverse Riemannian metric. The leaves correpond bijectively to the orbits of the holonomy pseudogroup $P$. If $M$ is a complete Riemannian manifold then the pseudogroup is complete as defined in [4] or [8], p. 278.

 An orthonormal coframe at $x \in X$ is an isometry $F_x$ from the tangent space $T_x$ of $X$ at $x$ to $R^n$. The coframes at various points of $X$ form a $O(n)$-principal bundle $p:F \to X$ over $X$. The space $F$ is a  Riemannian manifold  whose  tangent space  is parallelizable.
 The pseudogroup $P$ acts on $F$; this action defines a pseudogroup $P_F$ of local isometries of $F$ preserving the parallelism.  If $P$ is complete, then $P_F$ is also complete.
 
 For a pseudogroup like $P_F$ preserving a parallelism, the space of orbit closures is a manifold $W$. The space of orbits closures of $P$ is the quotient of $W$ by the action of $O(n)$. This space has a stratification according to the type of the orbits of $O(n)$ (see 
 M. Davis,
 Smooth $G$-manifolds as collections of fiber bundles,
  Pacific Journal of Math.  77,
( 1978), pp. 315--363).

 Also see [2] for a classification of the holonomy pseudogroups of leaf closures. Each case considered in this paper should be analyzed as in the present paper. 
 
 In a first version of this paper we used the notions of non-abelian cohomology to prove some results in section II,4. It would be an interesting exercise to generalize the concepts developped in the paper of Frenkel: Cohomologie non-ab\'elienne et espaces fibr\'es, Bull. Soc. Math. France, 8 (1957), pp 135-220,  by replacing, everywhere in this work, groups by groupoids, and also by introducing base points.


\bye

 \proclaim{Theorem} If the center $\Gamma_0$ of  $G$ is a discrete subgroup of $G$, then the classifying space  of   $(\Gamma \ltimes G- extension)$ is the classifying space of the Lie group $ Aut(\rho) \ltimes G /\Gamma_0$.

\endproclaim

\bye